\documentclass[10pt]{article}
\usepackage{amsfonts,a4wide}
\usepackage{url}
\usepackage[normalem]{ulem}
\usepackage{mathdots,amsmath,amssymb,amsthm,color}
\newtheorem{Th}{Theorem}[section]
\newtheorem{Lemma}[Th]{Lemma}

\newtheorem{proposition}[Th]{Proposition}
\newtheorem{Cor}[Th]{Corollary}

\newtheorem{Def}[Th]{Definition}

\theoremstyle{remark}
\newtheorem{Ex}{Example}
\newtheorem{Rem}[Th]{Remark}

\newcommand{\be}{\begin{equation}}
\newcommand{\ee}{\end{equation}}

\newcommand{\tr}{\mathrm{tr}\,}

\newcommand{\R}{\mathbb{R}}
\newcommand{\Id}{\textrm{\rm \bf 1}}

\newcommand{\Image}{\mathrm{Im}\,}
\newcommand{\trace}{\mathrm{tr}\,}

\newcommand{\weg}[1]{}

\title{Local normal forms for geodesically equivalent pseudo-Riemannian metrics}
\author{Alexey V. Bolsinov\footnote{ School of Mathematics,
 Loughborough University,
 LE11 3TU, UK  \ \
 \quad {\tt A.Bolsinov@lboro.ac.uk} } \footnote{Partially supported by Ministry of Education and Science of the Russian Federation (14.B37.21.1935)}  \qquad
\& \qquad  Vladimir S. Matveev\footnote{
Institute of Mathematics,
07737,  Jena,  Germany  \ \ \quad {\tt  vladimir.s.matveev@gmail.com}} \footnote{Partially supported by DFG (GK 1523) and DAAD (Programm Ostpartnerschaft) }}
\date{}
\begin{document}
\maketitle

\begin{abstract}
Two pseudo-Riemannian  metrics $g $ and $\bar g$ are geodesically equivalent, if they share  the same (unparameterized) geodesics.  We give a complete local description of  such metrics  which solves the natural generalisation of the Beltrami problem for  pseudo-Riemannian  metrics.
 \end{abstract}


\section{Introduction}


\subsection{Definition and  history} \label{sec:definitions}

Two pseudo-Riemannian metrics $g$ and $\bar g$
 on one manifold $M^n$ are geodesically equivalent, if every $g$-geodesic, after an appropriate reparameterisation, is a $\bar g$-geodesic. The theory of geodesically equivalent metrics had a long and rich history, the first examples being constructed by Lagrange in 1779 \cite{Lagrange}. Many important results about geodesically equivalent metrics were obtained by  Beltrami \cite{Beltrami, Beltrami2, Beltrami3}, Levi-Civita \cite{Levi-Civita}, Painlev\'e \cite{Painleve}, Lie \cite{Lie},  Liouville \cite{Liouville}, Fubini \cite{Fubini},   Eisenhart \cite{Eisenhart1,Eisenhart2}, Weyl \cite{Weyl} and  Thomas and Veblen \cite{thomas, veblen, thomasveblen}.
 Between 1950 and 1990, the theory of geodesically equivalent metrics was one of the main research areas of the Soviet and Japanese differential geometry schools, see the surveys \cite{Am2,mikes}. Recently, the theory of geodesically equivalent metrics has had a revival, in particular,  because of new mathematical
 methods that came from the theory of integrable systems \cite{MT}  and parabolic Cartan geometry \cite{BDE,eastwood}. Using these methods has led, in the last ten years, to the solution of many classical problems  including the Lie problems \cite{bryant,alone}, the Lichnerowicz  conjecture \cite{archive} and the Weyl-Ehlers problems \cite{einstein,relativity}.

 In this paper we solve the natural generalization of  a   problem explicitly stated by  Beltrami to the case of pseudo-Riemannian metrics, namely

  {\bf  Beltrami Problem\footnote{ Italian original from \cite{Beltrami}: \emph{
La seconda $\dots$  generalizzazione $\dots$ del nostro problema,   vale a dire:   riportare i punti di una superficie sopra un'altra superficie in modo  che alle linee geodetiche della prima corrispondano linee geodetiche della seconda}}.} \emph{ Describe all pairs of geodesically equivalent metrics.}

From the context it is clear that Beltrami actually considered this problem locally and in a neighborhood of almost every point, so do we. From the context it is also clear that   Beltrami thought about two dimensional  Riemannian surfaces; our answer does not have this restriction: the dimension of the manifold and the signatures of the metrics are arbitrary.

Special cases of the Beltrami  problem have previously been solved. The two-dimensional Riemannian case was solved by Dini \cite{Dini} in 1869. He has shown that two geodesically equivalent Riemannian metrics on a surface in  a neighborhood of almost every point are given,  in a certain coordinate system,
by the following formulas
\begin{equation}
\label{dini}
  g= (Y(y) - X(x))(dx^2 +dy^2)   \quad \textrm{and} \quad
\bar g=    \left(\frac{1}{X(x)}  - \frac{1}{Y(y)}\right)
    \left(\frac{dx^2}{X(x)} +\frac{dy^2}{Y(y)}\right).
\end{equation}
    Here $X$ and $Y$ are   functions of the indicated variables. For every smooth
    functions $X,Y$
     such that the formulas \eqref{dini} correspond to Riemannian metrics
     (i.e., $ 0<X(x)< Y(y)$  for all $(x,y)$), the metrics $g$ and $\bar g$ are geodesically equivalent.

For an arbitrary dimension, the Beltrami problem in the Riemannian case was solved by Levi-Civita \cite{Levi-Civita}. We will recall the Levi-Civita's  3-dimensional analog of the formulas \eqref{dini} below, in Example \ref{ex:dim3},  page \pageref{ex:dim3}. \weg{
Formally speaking, Levi-Civita solved only special cases of the Riemannian Beltrami problem: he gave the answer under the additional assumption that (1,1)-tensor $G$ given by the condition $g(G\,\cdot\, , \cdot)=\bar g(\cdot\, , \cdot)$
has $n=dim(M)$ different eigenvalues, or   $k-1$ different eigenvalues of multiplicity $1$ and $1$  eigenvaue of multiplicity $n-k$. The methods of Levi-Civita do not require this assumption and work for arbitrary structure of the tensor $G$ without essential modification. }

The methods of Levi-Civita and Dini can not be directly  generalized   to the pseudo-Riemannian case. Levi-Civita and Dini consider the tensor $G$ defined by the condition $g(G\,\cdot\, , \cdot)=\bar g(\cdot\, , \cdot)$. Levi-Civita has shown, that the eigenspaces of this tensor are simultaneously integrable which implies that (in a neighborhood of almost every point)  there exists a local coordinate system
$$
(\bar x_1, \dots ,\bar x_n)= ({x_1^1, \dots , x_{1}^{m_1}}, \, \dots \, ,x_{k}^{1},\dots ,x_{k}^{m_k})
$$
 such that  in these coordinates the matrix of $g$ is blockdiagonal with $k$ blocks of dimension $m_1, m_2,\dots ,m_k$ and the matrix of  $G$ is diagonal $\operatorname{diag}(\underbrace{\rho_1,\dots,\rho_1}_{m_1}, \, \dots \, ,\underbrace{\rho_k,\dots ,\rho_k}_{m_k})$.  In the two-dimensional Riemannian case, considered by Dini, the existence of such coordinate system is obvious.
 Now, in this  coordinate system, the partial differential equations  on the entries of $g$ and on  $\rho_i$ expressing the geodesic equivalence condition for  $g$ and $\bar g= g(G\,\cdot\, , \cdot)$ are relatively easy (though they are still coupled)  and, after some  nontrivial work, can be solved.

The methods of Levi-Civita also work in the pseudo-Riemannian  case under the additional assumption that $G$  is diagonalizable. Unfortunately, in the pseudo-Riemannian case the tensor $G$ may have complex eigenvalues and nontrivial Jordan blocks.
From the point of view of partial differential equations,
the case  of many Jordan blocks poses the main difficulties: unlike the case when $G$ is diagonalizable,  there is no `best' coordinate system, and  the equation corresponding to  the entries of the metrics coming from
different   blocks are coupled in a very nasty manner.

This difficulty was overcome in \cite{splitting}. In \S \ref{sec:splitting} we recall the main result of \cite{splitting} and  explain that
 the description of geodesically equivalent metrics $g$ and $\bar g$ in a neighborhood of almost every point can be reduced to the case when
the tensor $G$ has only one real eigenvalue, or two complex conjugate
eigenvalues.   The biggest part of our paper is devoted to the
 local description of geodesically equivalent metrics  under this assumption.

Special cases of the local description of geodesically equivalent pseudo-Riemannian metrics  were known before. The 2-dimensional case was described essentially by Darboux \cite[\S\S  593,  594]{Darboux}, see also \cite{pucacco,appendix}. Three dimensional case was solved by Petrov \cite{Petrov}, it is one of the results for which Petrov  obtained the Lenin prize  in 1972,
       the most important  scientific award of the  Soviet Union.  According to \cite{Am2}, under the additional assumption that the metrics $g$ and $\bar  g$ have Lorentz signature, the Beltrami problem was solved by
Golikov \cite{Golikov} in dimension 4, and by Kruchkovich \cite{Kruchkovich} in all dimensions; unfortunately,  we were not  able to find and to verify  these references.

It was generally believed that  the Beltrami problem was  solved in full generality in \cite{aminova}. Unfortunately, this  result of Aminova  seems to be wrong.   More precisely,   in view of \cite[Theorem 1.1]{aminova}   and  the formulas  \cite[(1.17),(1.18)]{aminova} for $k=1$, $n=4$ and all $\varepsilon$'s equal to $+1$,
   the following two metrics $g$ and $\bar g$ given by the matrices  (where $\omega$ is an arbitrary function of the variable $x_4$)
$$
\left[ \begin {array}{cccc}
0&0&0&3\,x_{{3}}+3\,\omega \left( x_{{4}}
 \right) \\\noalign{\medskip}0&0&1&2\,x_{{2}}\\\noalign{\medskip}0&1&0
&x_{{1}}\\\noalign{\medskip}3\,x_{{3}}+3\,\omega \left( x_{{4}}
 \right) &2\,x_{{2}}&x_{{1}}&4\,x_{{1}}x_{{2}}
 \end {array}
 \right] ,
 $$

 {\tiny $
 \left[ \begin {array}{cccc} 0&0&0&3\,{\frac {x_{{3}}+\omega \left( x_
{{4}} \right) }{{x_{{4}}}^{5}}}\\\noalign{\medskip}0&0&2\,{x_{{4}}}^{-
5}&{\frac {-3\,x_{{3}}-3\,\omega \left( x_{{4}} \right) +2\,x_{{2}}x_{
{4}}}{{x_{{4}}}^{6}}}\\\noalign{\medskip}0&2\,{x_{{4}}}^{-5}&-{x_{{4}}
}^{-6}&{\frac {3\,x_{{3}}+3\,\omega \left( x_{{4}} \right) -2\,x_{{2}}
x_{{4}}+x_{{1}}{x_{{4}}}^{2}}{{x_{{4}}}^{7}}}\\\noalign{\medskip}3\,{
\frac {x_{{3}}+\omega \left( x_{{4}} \right) }{{x_{{4}}}^{5}}}&{\frac
{-3\,x_{{3}}-3\,\omega \left( x_{{4}} \right) +2\,x_{{2}}x_{{4}}}{{x_{
{4}}}^{6}}}&{\frac {3\,x_{{3}}+3\,\omega \left( x_{{4}} \right) -2\,x_
{{2}}x_{{4}}+x_{{1}}{x_{{4}}}^{2}}{{x_{{4}}}^{7}}}&{\frac { \left( -3
\,x_{{3}}-3\,\omega \left( x_{{4}} \right) +2\,x_{{2}}x_{{4}} \right)
 \left( 2\,x_{{1}}{x_{{4}}}^{2}+3\,x_{{3}}+3\,\omega \left( x_{{4}}
 \right) -2\,x_{{2}}x_{{4}} \right) }{{x_{{4}}}^{8}}}\end {array}
 \right]$}

should be geodesically equivalent, though they are not (which can be checked by direct calculations).



\subsection{Splitting and gluing construction: why it is sufficient to assume that  {\mathversion{bold} $G$} has one real eigenvalue  or two complex conjugate eigenvalues.} 
\label{sec:splitting}


Given two metrics $g$ and $\bar g$ on the same manifold, instead of considering the $(1,1)$-tensor
$G^i_j= g^{ik}\bar g_{kj}$,
 we consider the $(1,1)-$tensor $L=L(g,\bar g)$ defined by
  
\begin{equation}
\label{L}
L_j^i := { \left|\frac{\det(\bar g)}{\det(g)}\right|^{\frac{1}{n+1}}} \bar g^{ik}
 g_{kj},
\end{equation}
where ${\bar g}^{ik}$ is the contravariant inverse of
${\bar g}_{ik}$. The tensors $G$ and $L$ are related by 
\begin{equation}
\label{LG}
L= { |\det(G)|}^\frac{1}{n+1}G^{-1},  \quad  G= \tfrac{1}{{ |\det(L)|}} L^{-1},
\end{equation}
 so  in particular they have the same structure
 of  Jordan blocks (though their eigenvalues are, in general, different). Since the metric $\bar g$  can be  uniquely reconstructed from $g$ and $L$, namely:
 \begin{equation}
 \label{gandL}
 \bar g (\cdot\, , \cdot) = \tfrac{1}{{ |\det(L)|}}g(L^{-1}\cdot\, ,\cdot)  
 \end{equation}
 the condition that $\bar g$  is geodesically equivalent to $g$ can be written as a system of PDEs on the components of $L$.     From   the point of view of partial differential equations, the tensor $L$ is more convenient than $G$: the corresponding system of
   partial differential equations on $L$ turns out to be linear. In the index-free form, it  can be written as the condition  (where ``$\ast$'' means  $g-$adjoint)
\begin{equation}
\label{main}
\nabla_u L =\frac{1}{2} (u\otimes d\tr L + (u\otimes d\tr L )^* ),
\end{equation}
which should be fulfilled at every point and for every vector field $u$.

In tensor notation, condition \eqref{main} reads
\begin{equation}
\label{main1}
L_{ij,k} = \lambda_{,i} g_{jk} + \lambda_{,j} g_{ik},
\end{equation}
where  $L_{ij}:= L^k_j g_{ki}$ and  $\lambda:= \tfrac{1}{2} L^i_i= \tfrac{1}{2} \trace(L)$.
The tensor $L^i_j $  defined  by  \eqref{L} is essentially the  same
 as the tensor   introduced by Sinjukov (see equations (32, 34) on page 134 of the book \cite{sinjukov} and also
Theorem 4 on page 135); the equation  \eqref{main1} is  also due to him, see also \cite[Theorem 2]{benenti}.

\weg{
\begin{Rem}
If $n$ is even, the tensor $L$ is always well defined. If $n$ is odd,  the ratio ${\det(\bar g)}/{\det(g)}$ may be negative and  then formula \eqref{L}  makes no sense. There  is  the following  way to avoid this (rather formal) difficulty:
 we  can
replace   $\bar g$  by $-\bar g$ and  make the ratio ${\det(\bar g)}/{\det(g)}$ positive  and $L$  well defined.     Moreover, since the equations \eqref{main} are linear, we can assume  that $L$ is close to $\Id$ in the neighborhood we are working in, implying
 that  $g$ and $\bar g$ have the same signatures, and the problem with the sign does not appear at all.
\end{Rem}}

\begin{Def}{\rm
\label{compatibility}
We say that a (1,1)-tensor  $L$ is {\it compatible} with $g$, if  $L$ is $g$-selfadjoint, nondegenerate at every point and satisfies \eqref{main} at any point and for all tangent vectors $u$.  
}\end{Def}

As we explained above, $L$ is compatible with $g$ if and only if $\bar g(\cdot\, , \cdot)= \tfrac{1}{{ |\det(L)|}} g(L^{-1}\cdot\, , \cdot)$ is a pseudo-Riemannian metric geodesically equivalent  to $g$.

 The {\it gluing construction},  as well as the {\it splitting construction} to be presented below,  are due  to \cite{splitting};
in the  Riemannian  case they appeared slightly earlier, see \cite[\S 4]{hyperbolic}, 
  \cite[Lemma 2]{archive} and  \cite[\S\S2.2, 2.3]{bifurcations}.

Consider two pseudo-Riemannian manifolds $(M_1, h_1)$ and $(M_2,h_2)$. Assume that $L_1$ on $M_1$  is  compatible with $h_1$,  and that  $L_2$ on $M_2$  is  compatible with $h_2$. Assume in addition that  $L_1$ and $L_2$
 have no common eigenvalues in the sense that
for  any two points $x\in M_1$, $y\in M_2$ we have
\begin{equation}
\label{spectrum}
\textrm{Spectrum}\, L_1(x) \cap \textrm{Spectrum}\, L_2(y) =\varnothing.
\end{equation}

Then one can naturally and canonically
construct a  pseudo-Riemannian metric $g$  and a tensor $L$ compatible with $g$ on the direct product   $M=M_1 \times  M_2$.
The new metric   $g$ differs  from the direct product metric  $h_1 + h_2$  on $M_1\times M_2$ and  is  given  by the following   formula involving $L_1$ and $L_2$.
We denote by $\chi_i$, $i=1,2$, the characteristic polynomial of $L_i$: $\chi_i= \det(t\cdot \Id_i - L_i)$ (where $\Id_i$ is the identity operator $\Id_i:TM_i\to TM_i$).   We treat the $(1,1)-$tensors $L_i$ as linear operators acting on $TM_i$.
  For a polynomial $f(t) =  a_0  + a_1 t + a_2 t^2  + \cdots + a_m  t^m$ and $(1,1)$-tensor $A$, we put   $f(A)$ to be  the $(1,1)$-tensor
  $$
  f(A)=a_0 \cdot \Id + a_1 A + a_2 A\circ A + \cdots +
  a_m \underbrace{A\circ\cdots\circ A}_{\textrm{\tiny $m $ times }}.
  $$

If no eigenvalue of $A$ is a root  of $f$, $f(A)$ is nondegenerate; if $A$ is $g$-selfadjoint, $f(A)$ is  $g$-selfadjoint as well.

For two tangent vectors
$
u= (\underbrace{u_1}_{\in TM_1}, \underbrace{u_2}_{\in TM_2})\, , \  \
v=( \underbrace{v_1}_{\in TM_1}, \underbrace{v_2}_{\in TM_2}) \in TM
$
 we put    \begin{eqnarray} g(u,v) &  = &  h_1\left( \chi_2(L_1)( u_1), v_1\right)   + h_2\left(\chi_1(L_2)(u_2), v_2\right),  \label{hh1}  \\
  L(v) &= & \left(L_1(v_1), L_2(v_2)\right).  \label{Lnew}
\end{eqnarray}
We see that   the   $(1,1)-$tensor $L$ is  the direct sum of $L_1$ and $L_2$  in the natural sense.

It might   be  convenient to understand  the formulas (\ref{hh1}, \ref{Lnew})  in  matrix notation: we  consider   the coordinate system
$(x_1,\dots ,x_r,y_{r+1},\dots ,y_{n})$ on $M$  such that  $x$'s are coordinates on    $M_1$ and  $y$'s are coordinates on  $M_2$.
Then, in this coordinate system, the matrices  of $g$ and $L$   have the block diagonal form

\begin{equation}
\label{matg}
g =\begin{pmatrix}   h_1   \chi_2(L_1) & 0 \\  0 &  h_2   \chi_1(L_2)\end{pmatrix}\ , \ \ L =\begin{pmatrix} L_1 & 0 \\  0 &  L_2\end{pmatrix}.
\end{equation}

If  \eqref{spectrum} is fulfilled, then $g$ is a pseudo-Riemannian metric (i.e., symmetric and nondegenerate) and   $L$  is nondegenerate and  $g$-selfadjoint.

\begin{Th}
[Gluing Lemma from \cite{splitting}] \label{thm3}
If $L_1$ is compatible with $h_1$ on $M_1$,   $L_2$ is compatible with   $h_2$ on $M_2$ and  \eqref{spectrum} is  fulfilled,   then  $L$  given by \eqref{Lnew}  is compatible with $g$ given by \eqref{hh1}.
\end{Th}

\begin{Ex}
[Gluing construction and Dini's Theorem]
As the manifolds $M_1,M_2$ we take two intervals $I_1$ and $I_2$ with coordinates $x$ and  $y$ respectively. Next, take the metrics $h_1= dx^2$ on $I_1$ and $h_2= - dy^2$ on $I_2$. Consider the $(1,1)$-tensors
 $L_1= X(x) dx\otimes \tfrac{\partial}{\partial x} $ on $I_1$ and
 $L_2= Y(y) dy\otimes \tfrac{\partial}{\partial y} $ on $I_2$. We assume that  $0<X(x)<Y(y)$ for all $x\in I_1$ and $y\in I_2$.   The  tensors $L_1$ and $L_2$  are compatible  with $h_1$ and $h_2$ respectively (which can be checked by direct calculation and which is trivial in view of the obvious fact that in dimension $1$ all metrics are geodesically equivalent). Then, $\chi_{1}= (t- X(x))$,   $\chi_{2}= (t- Y(y))$, so the formulas \eqref{matg} read
 $$
 g =
 \begin{pmatrix}
 Y(y)-X(x) & \\ & Y(y)-X(x)
 \end{pmatrix}, \
 L =
 \begin{pmatrix}
 X(x)& \\ &Y(y)
 \end{pmatrix}.
 $$
 Now, combining this with \eqref{gandL}, we obtain that this  $g$  is geodesically equivalent to the metric
 $$
 \bar g =\frac{Y(y)-X(x)}{X(x)Y(y)}
 \begin{pmatrix}
 \tfrac{1}{X(x)}& \\ & \tfrac{1}{Y(y)}
 \end{pmatrix}
 = \left(\frac{1}{X(x)}- \frac{1}{Y(y)}\right)
 \begin{pmatrix}
 \tfrac{1}{X(x)}& \\ & \tfrac{1}{Y(y)}
 \end{pmatrix}.
 $$
 Comparing the above formulas with \eqref{dini}, we see that  the gluing construction  applied
 to two intervals proves  Dini's local description of geodesically equivalent metrics in one direction.
 \end{Ex}

One can iterate this construction: having three pseudo-Riemannian manifolds $(M_1, h_1)$, $(M_2, h_2)$, $(M_3,h_3)$ carrying $h_i$-compatible $(1,1)$-tensors $L_i$ with pairwise disjoint spectra (see  \eqref{spectrum}), one can glue $M_1$ and $M_2$ and then glue the result with $M_3$. Actually the  gluing construction is  associative. Indeed,  one obtains the same metric   $g$ and the same $g$-compatible $(1,1)$-tensor $L$ on $M_1+M_2+M_3$ if one
 first glues $(M_1, h_1, L_1)$  and $(M_2, h_2, L_2)$ and then glues  the result with $(M_3,h_3, L_3)$,
  or if one first glues $(M_2, h_2 L_2)$ and  $(M_3,h_3, L_3)$ and then glues  $(M_1,h_1, L_1)$ with the result:
  $$
  \bigl((M_1,h_1, L_1) \stackrel{\textrm{\tiny glue}}{+}(M_2, h_2, L_2) \bigr) \stackrel{\textrm{\tiny glue}}{+} (M_3,h_3, L_3) = (M_1,h_1, L_1) \stackrel{\textrm{\tiny glue}}{+}  \bigl((M_2,h_2, L_2) \stackrel{\textrm{\tiny glue}}{+}(M_3, h_3, L_3) \bigr).
  $$
The gluing construction is commutative as well:
$$
(M_1,h_1, L_1) \stackrel{\textrm{\tiny glue}}{+}(M_2, h_2, L_2) \stackrel{\textrm{\tiny iso}}{=} (M_2,h_2, L_2) \stackrel{\textrm{\tiny glue}}{+}(M_1, h_1, L_1),
$$
where
``$\stackrel{\textrm{\tiny iso}}{=}$'' means the existence of a diffeomorphism that preserves the metric and $L$. Actually, this diffeomorphism is given by the  natural formula
$$
M_1\times M_2 \ni (\underbrace{x}_{\in M_1},\underbrace{y}_{\in M_2})
\mapsto
(\underbrace{y}_{\in M_2},\underbrace{x}_{\in M_1}) \in M_2\times M_1.
$$

In the case we ``glue'' $k$ manifolds   $(M_i,h_i)$ ($i= 1,\dots ,k$) such that each manifold is equipped with $h_i$-compatible $L_i$, we obtain a metric $g$  on $M= M_1\times \cdots \times M_k$ and $g$-compatible $L$ on  $M$ such that in the matrix notation in the natural coordinate system  they  have the form
\begin{equation}
\label{manyblocksg}
 g=
 \textrm{$
 \begin{pmatrix}
 h_1 \chi_2(L_1)\weg{\chi_3(L_1)}\cdots\chi_k(L_1) & & & \\
                      & h_2  \chi_1(L_2)\chi_3(L_2)\cdots\chi_k(L_2)&& \\
                      & & \ddots &\\
                      & & & h_k  \chi_1(L_k)\weg{\chi_2(L_k)}\cdots\chi_{k-1}(L_k)
\end{pmatrix}
$}
\end{equation}
and
\begin{equation}
\label{manyblocksL}
L=
\begin{pmatrix}
L_1 & & & \\
        & L_2&& \\
        & & \ddots &\\
        & & & L_k
\end{pmatrix}.
\end{equation}

\begin{Ex}
[Gluing construction and Levi-Civita's Theorem in dim 3]
\label{ex:dim3}
We now take three intervals
  $I_1, I_2, I_3$ with the coordinates $x$, resp. $y$, $z$,    the metrics $h_1= dx^2$, $h_2= - dy^2$, $h_3= dz^3$,  and the $h_i$-compatible $(1,1)$-tensors
 $L_1= X(x) dx\otimes \tfrac{\partial}{\partial x} $,
 $L_2= Y(y) dy\otimes \tfrac{\partial}{\partial y} $ and $L_3= Z(z) dz\otimes \tfrac{\partial}{\partial z} $. We  again assume that  the spectra  of $L_i$ are pairwise disjoint at every point and 
 $$
0<X(x)<Y(y)<Z(z) \quad \forall x\in I_1, \ \forall y\in I_2, \ \forall z\in I_3.
 $$
 Applying the gluing construction twice, we obtain the metric 
 $g$ and the $(1,1)$-tensor   $L$ on $M^3= I_1\times I_2 \times I_3$  that are compatible  and given in the natural coordinate system by:
$$
 g =
 \textrm{$
 \begin{pmatrix}
 (Y(y)\!-\!X(x))(Z(z)\!-\!X(x)) & & \\
 & \!\!\! (Y(y)\!-\!X(x))(Z(z)\!-\!Y(y))& \\
 &&  \!\!\! (Z(z)\!-\!Y(y))(Z(z)\!-\!X(x)) \\
 \end{pmatrix}$} 
 $$
 and
 $$
 L=
 \textrm{$
 \begin{pmatrix}
 X(x)& &\\
 &Y(y)&\\
 &&Z(z)
 \end{pmatrix}
 $}.
 $$

  Combining this with   \eqref{gandL}, we obtain a special case of  Levi-Civita's  local description of geodesically equivalent metrics in dimension 3 (when the tensor $G$ has three different eigenvalues).
 \end{Ex}

The {\it splitting construction} is the inverse operation.  We will   describe its  local version only, since it is sufficient for our goals.

Suppose $g$ is a pseudo-Riemannian metric on $M^n$ and $L$ is compatible with $g$.  We consider an arbitrary point $p\in M$.

We take a point $p$ of the manifold such that in the neighborhood $U(p)$  of this point the eigenvalues of $L$ do not bifurcate (i.e., the number of different eigenvalues is constant in the neighborhood). Then, the eigenvalues are smooth possibly complex-valued functions. We denote them by
$$
\lambda_1,\bar \lambda_1, \, \dots \, ,\lambda_r,\bar \lambda_r: U(p)\to \mathbb{C}, \quad
 \lambda_{r+1},\, \dots \, ,\lambda_k:U(p)\to \mathbb{R}.
 $$
We assume that  for $i\le r$  the eigenvalue
$\bar \lambda_i$ is complex conjugate to $\lambda_i$ and $\mathrm{Im}\, \lambda_i\ne 0$. We think that the eigenvalue $\lambda_i$ has algebraic multiplicity $m_i$, $2m_1+ \cdots +2m_r + m_{r+1}+\cdots + m_{k}= n$.

 Next, let us consider the polynomial functions $\chi_i:\mathbb{R}\times U(p) \to \mathbb{R}$:
 $$
 \chi_i = (t- \lambda_i)^{m_i}(t-\bar \lambda_i)^{m_i} \ \  \textrm{for $i=1,\dots ,r$ \  and } \ 
 \chi_i = (t- \lambda_i)^{m_i} \ \  \textrm{for $i=r+1,\dots \dots ,k$, }
 $$
 and the polynomial function $\hat \chi:= \hat\chi_1+\dots+\hat\chi_k$, where $\hat\chi_i = \frac{\chi}{\chi_i}$ and $\chi= \det(t\cdot \Id - L)$ is the characteristic polynomial of $L$. It is easy to see  that the $(1,1)$-tensor
 $\hat\chi(L)$ is $g$-selfadjoint and  nondegenerate. Then we can introduce a new pseudo-Riemannian metric $h$ on $U(p)$ by
\begin{equation}\label{222}
 h(u,v):=  g(\hat\chi(L)^{-1}u,v),  \qquad u,v\in T_qM, \ \  q\in U(p).
 \end{equation}

  \begin{Th}
  [Follows from the Splitting Lemma, see \S 2.1 of \cite{splitting}]
  \label{thm:splitting}
  In a neighborhood of $p$ there exists a coordinate system
$$
 (\bar x_1,\dots ,\bar x_k)= \bigl(x_1^1, \dots ,x_1^{2m_1},\cdots,x_r^1,\dots ,
 x_r^{2m_r},x_{r+1}^1,\dots ,x_{r+1}^{m_{r+1}},\cdots,
 x_{k}^1,\dots ,x_{k}^{m_{k}}\bigr)
$$
 in which the matrices of $h$ and of $L$ are given by
 \begin{equation}
 \label{manyblocksLtheorem}
 h=  \begin{pmatrix} h_1 & & & \\
                      & h_2&& \\
                       && \ddots &\\
                         &&& h_k
\end{pmatrix}, \ \
L=  \begin{pmatrix} L_1 & & & \\
                      & L_2&& \\
                       && \ddots &\\
                         &&& L_k
\end{pmatrix}.
\end{equation}
  Moreover,
  \begin{itemize}
  \item the entries of the blocks $h_i$ and $L_i$ depend on the coordinates $\bar x_i$ only;

  \item for $i=1,\dots ,r$ the eigenvalues of $L_i$ are $\lambda_i$ and $\bar \lambda_i$,
  and for  $i=r+1,\dots ,k$ the only  eigenvalue of $L_i$ is  $\lambda_i$;

   \item $L_i$ is compatible with $h_i$ for every $i=1,\dots ,k$. 
    \end{itemize}
  \end{Th}

 \begin{Ex}
 [Splitting  construction and Dini's Theorem]
 Consider a two dimensional manifold $M^2$   and geodesically equivalent Riemannian metrics $g$ and $\bar g$ on it. Let $p\in M$ be a point where the metrics are not proportional. Then, $L(g,\bar g)$ has two (real)
 eigenvalues $\lambda_1 \not = \lambda_2$ at every point of a small neighborhood $U(p)$. Then,  in the notation of Theorem \ref{thm:splitting}, $k=2$, $r=0$ and $m_1=m_2=1$.  Thus, $\chi_1= t-\lambda_1$, $\chi_2= t-\lambda_2$ and $\hat \chi = (t- \lambda_2)+ (t- \lambda_1)$.    Then, there exists a coordinate system $x,y$
 such that $\lambda_1=X(x)$, $\lambda_2=Y(y)$ and  $h$, $L$ and $\hat\chi(L)$ are given by the  matrices
 $$
 h =
 \begin{pmatrix}
 \tilde X(x) & \\
 & \tilde Y(y)
 \end{pmatrix}, \  \
 L=
 \begin{pmatrix}
 X(x) & \\ & Y(y)
 \end{pmatrix}, \ \ 
 \hat\chi( L)=
 \begin{pmatrix}
 Y(y)- X(x) & \\ & X(x)-Y(y)
 \end{pmatrix}.
 $$
 Combining this with  \eqref{gandL}, \eqref{222}, we see that the metrics $g$ and $\bar g$ are given by
$$
 g= (Y(y) - X(x))(\tilde X(x) dx^2 + \tilde Y(y) dy^2)  \ \ \textrm{and} \ \
  \bar g =  \left(\frac{1}{X(x)}  - \frac{1}{Y(y)}\right)\left(\frac{\tilde X(x)dx^2}{X(x)} +\frac{\tilde Y(y) dy^2}{Y(y)}\right).
$$
By a coordinate change of the form $x_{\mathrm{new}}= x_{\mathrm{new}}(x), y_{\mathrm{new}}= y_{\mathrm{new}}(x)$ one can `hide' $\tilde X$ in $dx^2$ and $ \tilde Y$  in $dy^2$ and obtain the metrics of the form
 \eqref{dini}.
 \end{Ex}

{\it Vocabulary.} \label{vocal} We call a point $p \in M$ {\it regular},  if in some  neighborhood $U(p)$ of $p$ the Jordan  type  of $L$ is constant (that is, the number of eigenvalues and  Jordan  blocks is the same at all points $x\in U(p)$; the sizes of Jordan blocks are assumed to be fixed too, whereas the eigenvalues can, of course, depend on the point).
  It is easy to see that almost every point of $M$ is regular, that is,  the set of regular points is open and dense on the manifold.

  Now, applying the Splitting Lemma in the neighborhood of a regular point, we obtain the metrics $h_i$ on $2m_i$- or $m_i$  dimensional discs and $(1,1)$-tensors $L_i$ compatible with $h_i$. Moreover, each $L_i$ has one real or two complex conjugate eigenvalues and the Jordan type  of  $L_i$ is the same at all points.

  Describing all possible   $h_i$ and $L_i$ satisfying these conditions gives
   all possible geodesically equivalent  metrics $g$ and $\bar g$ near  regular points:
  all possible $g$ and $g$-compatible $L$ can be obtained  by the gluing construction (which is given by explicit formulas (\ref{manyblocksg},\ref{manyblocksL}) and the metric $\bar g$ is constructed from $g$ and $L$ by the formula  \eqref{gandL}.

  Finally,  in order to describe the  metric  and $L$ in the neighborhood of almost any point, it is sufficient to describe the metrics $h_i$ and the $h_i$-compatible $L_i$ such that  $L_i$  has one real  eigenvalue, or two complex conjugate eigenvalues, and the type of the Jordan block is the same in the whole neighborhood.  We will formulate the result in \S \  \ref{canonical_form}, see Theorems \ref{multtwoth}, \ref{onerealblock}, \ref{onecomplexblock} there, the  proof of these theorems will be given    in Sections \ref{affine}, \ref{proof2} and \ref{proof3}.



\subsection{Canonical forms for basic blocks} \label{canonical_form}


Throughout this section we assume that $p\in M$ is a regular point, i.e., the Jordan type of $L$ remains unchanged in some neighborhood of $p$. Our goal is to find local normal forms for compatible $L$ and $g$ nearby $p$ (see Definition~\ref{compatibility}).

According to the previous section (Theorem \ref{thm:splitting}), it is sufficient to describe the structure of compatible pairs $(g,L)$ in the case when $L$ either has a single real eigenvalue $\lambda$, or has a pair of complex non real eigenvalues $\lambda, \bar\lambda$. However even in these cases, the situation depends essentially  on    the algebraic structure of $L$, more precisely on  the geometric multiplicity of $\lambda$, i.e., the number of linearly independent eigenvectors. There are three essentially different possibilities: 1) the geometric multiplicity of  $\lambda$  is at least two,   
2) $L$ is conjugate to a real Jordan block (i.e., $L$ has a single real eigenvalue $\lambda$ of geometric multiplicity one) and 3)  $L$ is conjugate to a pair of complex conjugate Jordan blocks (i.e., $L$ has a pair of complex conjugate eigenvalues each of geometric multiplicity one). These cases are described by Theorems \ref{multtwoth},  \ref{onerealblock} and \ref{onecomplexblock} below.

We start with the case of multiplicity $\ge 2$. This situation turns out to be very special. Namely, the following statement holds.

\begin{Th}\label{multtwoth}
Let $g$ and $L$ be compatible and, in a neighborhood  $U$ of a point $p\in M$, the operator $L$ have either a unique real eigenvalue $\lambda=\lambda(x)$ or a unique pair of complex conjugate eigenvalues $\lambda(x), \bar\lambda(x)$, $\mathrm{Im}\, \lambda (x)\ne 0$. Suppose that the geometric multiplicity of $\lambda$ is at least two at each point $x\in U$.     Then  the function $\lambda (x)$ is constant and $L$ is covariantly constant in $U$, i.e., $\nabla L=0$. In particular the metrics $g$ and $\bar g$ given by \eqref{gandL} are affinely equivalent.
\end{Th}

Thus, in the case of geometric multiplicity $\ge 2$, our problem is reduced to another rather  nontrivial problem of local classification of pairs $g$, $L$
satisfying $\nabla L=0$, which has been recently completely solved by Charles Boubel and we refer to his work \cite{boubel}  for further details.

We now give the answer for $L$ being conjugate to a single real Jordan block, in other words we assume  that the eigenvalue $\lambda$ is real and $L$ possesses a unique (up to proportionality) eigenvector.

\begin{Th}\label{onerealblock}
Let $g$ and $L$ be compatible and $L$ be conjugate to a single Jordan block with a real eigenvalue $\lambda$. Then there exists a
local coordinate system $x_1,\dots, x_n$ such that $\lambda$ depends only on $x_n$ and that:
$$
g=\begin{pmatrix}
 & & &   & a_{n-1} \\
 & & & 1 & a_{n-2} \\
 & & \iddots &      & \!\!\!\!\!\!\!\!\!\!\!\!\!\! \!\!\!\!\!\!\!\!\!\! \mbox{{\huge $0$}} \ \ \quad\quad \vdots \\
 & 1 &     & & a_1 \\
a_{n-1} & a_{n-2} & \dots & a_1 & \sum^{n-2}_{i=1} a_i a_{n-i-1}
\end{pmatrix}
$$
and
$$
L=\begin{pmatrix}
\lambda(x_n) & 1            &       & &  a_1 \\
             & \lambda(x_n) & \ddots     & &  a_2 \\
             &              & \ddots &  1 &   \vdots  \\
             &              &        & \lambda(x_n) & a_{n-1}         \\
             &              &        &              & \lambda(x_n)
\end{pmatrix}
$$
where
$$
\aligned
& a_1=\lambda_{x_n}'  x_1,  \\
& a_2=2\lambda_{x_n}' x_2,  \\
& \dots                              \\
& a_{n-2}=(n-2)\lambda_{x_n}' x_{n-2}, \\
& a_{n-1}= 1 + (n-1)\lambda_{x_n}' x_{n-1},
\endaligned
$$
and $\lambda_{x_n}' $ stands for $\partial \lambda/\partial x_n$.
 Conversely, if $\lambda=\lambda(x_n)$ is an arbitrary smooth function such that $\lambda(x_n)\ne 0$ for all $x_n$, then $g$ and $L$ given by the above formulas are compatible (in the domain where $g$ is nondegenerate, i.e.  $1 + (n-1)\lambda_{x_n}' x_{n-1}\ne 0$).
\end{Th}

\begin{Rem}
Equivalently, one can write $g$ as the symmetric 2-form
$$
\sum_{k=1}^n \bigl(dx_{k} + (k-1)\, \lambda'_{x_n}  x_{k-1} dx_n\bigr)\bigl(dx_{n-k+1} + (n-k)\, \lambda'_{x_n}  x_{n-k} dx_n\bigr),
$$
with the convention $x_0=0$.  
The $(1,1)$-tensor  $L$, in this notation, takes the following form:
$$
L = \lambda(x_n)\cdot \Id + \sum_{k=1}^{n-1} \partial_{x_k}\otimes dx_{k+1} + \lambda'_{x_n}
\left( \sum_{k=1}^{n-1} kx_k \partial_{x_k} \right)\otimes dx_n.
$$

\end{Rem}

\begin{Rem}
In the case when $\lambda_{x_n}'\ne 0$ at the point $p$, we can simplify these formulas even further by taking $\lambda(x_n)$ as a new coordinate.
After the change $x^{\mathrm{new}}_n=\lambda(x_n)$ we obtain the following normal forms for $L$ and $g$ (we keep the ``old'' notation $x_n$ for the ``new'' coordinate).

Let $g$ and $L$ be compatible and $L$ be conjugate to a Jordan block with a real eigenvalue $\lambda$. If $d\lambda(p)\ne 0$, then in a neighborhood of $p\in M$ there exists a
local coordinate system $x_1,\dots, x_n$ such that  $\lambda = x_n$, $\lambda_0 = x_n(p)\ne 0$, and
\begin{equation}
\label{addg}
g=
\begin{pmatrix}
 & & &   & h(x_n)\!+\!(n\!-\!1)x_{n-1} \\
 & & & 1 & (n\!-\! 2) x_{n-2} \\
 & & \iddots & & \!\!\!\!\!\!\!\!\!\!\!\!\!\! \!\!\!\!\!\!\!\!\!\! \mbox{{\huge $0$}} \ \ \quad\quad \vdots \\
 & 1 &     & & x_1 \\
h(x_n)\!+\!(n\!-\!1)x_{n-1}& \!(n\!-\! 2) x_{n-2} & \dots & x_1 & \sum
\end{pmatrix}
\end{equation}
and
\begin{equation}
\label{addL}
L=\begin{pmatrix}
x_n & 1            &       & &  x_1 \\
             & x_n & \ddots     & &  2x_2 \\
             &              & \ddots & 1 &  \vdots  \\
             &              &        & x_n & h(x_n)\!+\!(n\!-\!1)x_{n-1}         \\
             &              &        &              & x_n
\end{pmatrix}
\end{equation}
where  $\sum = \sum^{n-2}_{i=1} i (n-i+1) x_i x_{n-i-1}$ and $h(x_n)$ is an arbitrary function such that $h(\lambda_0)\ne 0$ (equal to $1/\lambda'_{x^{\mathrm{old}}_n}$). Conversely, $g$ and $L$ given by these formulas are compatible for every $h(x_n)$  (in the domain where $x_n\ne 0$ and $h(x_n)+(n-1)x_{n-1}\ne 0$).
\end{Rem}

\begin{Rem} It follows immediately from the proof (see Section \ref{proof2}) that the canonical coordinate system (and hence canonical forms) for $g$ and $L$ from Theorem \ref{onerealblock}   is uniquely defined (up to a finite group) if we fix the position of the initial point $p\in M$ by saying that coordinates of $p$ in the canonical coordinate system are $(0,\dots, 0, \lambda_0)$.   If we do not fix $p$, i.e., move the origin to another point $p'\in U(p)$, then the function $h(x_n)$ playing the role of the parameter for canonical forms \eqref{addg} and \eqref{addL}  changes.  It is not difficult to check that the transformation that preserves the structure of  \eqref{addg} and \eqref{addL} has the following form:
$$
\begin{array}{l}
\tilde x_n = x_n, \\
\tilde x_{n-1}=x_{n-1} + P (x_n), \\
\tilde x_{n-2}=x_{n-2} - \frac{1}{n-2} P'(x_n),\\
\dots \\
\tilde x_{n-k}= x_{n-k} + (-1)^{k-1} \frac{1}{(n-2)(n-3)\dots(n-k)} P^{(k-1)} (x_n),\\
\dots \\
\tilde x_1 = x_1 + (-1)^{n-2}\frac{1}{(n-2)!} P^{(n-2)}(x_n).
\end{array}
$$
Here $P(x_n)$ is an arbitrary polynomial of degree $n-2$ and $P^{(k)}$ denotes its $k$th derivative. The function $h(x_n)$ after this change of variables takes the form $h(\tilde x_n)- (n-1) P(\tilde x_n)$. Thus we see that  the function $h$, the parameter of our canonical form, is defined modulo a polynomial of degree $n-2$.   
\end{Rem}

The next case is a complex Jordan block, i.e. we assume that the only eigenvalues of $L$ are a pair of complex conjugate numbers $\lambda$ and $\bar\lambda$ for each of which there is a single (up to proportionality) eigenvector over $\mathbb C$. Equivalently, this means that the corank of the real operator $(L-\lambda\cdot\Id)(L-\bar\lambda\cdot\Id)$ is two. In this case, the normal form for $g$ and $L$ can be described in the following way.

\begin{Th}\label{onecomplexblock}
Let $g$ and $L$ be compatible and $L$ be conjugate to a complex Jordan block with complex conjugate eigenvalues $\lambda$ and $\bar\lambda$  {\rm (}$\mathrm{Im}\, \lambda\ne 0${\rm )}. Then there exists a complex structure $J$ and a local complex coordinate system
$(z_1,\dots, z_n)$ such that
\begin{enumerate}
\item
the eigenvalue $\lambda$ is a holomorphic function of $z_n$,
\item
$L$ is a complex linear operator on $(T_PM, J)$ given in this coordinate system by the matrix:
$$
L^{\mathbb C} =
\begin{pmatrix}
\lambda (z_n) & 1            &        & &  a_1 \\
             & \lambda (z_n) & \ddots      & &  a_2 \\
             &               & \ddots &  1 &   \vdots  \\ 
             &               &        & \lambda(z_n) & a_{n-1} \\
             &               &        &              & \lambda (z_n)
\end{pmatrix}
$$
\item
the metric $g$ is the real part of the complex bilinear form on $(T_PM, J)$ given in this coordinate system by the matrix:
$$
g^{\mathbb C}=-i \begin{pmatrix}
& & & & a_{n-1} \\
 & & & 1& a_{n-2}\\
 & & \iddots & & \!\!\!\!\!\!\!\!\!\!\!\!\!\! \!\!\!\!\!\!\!\!\!\! \mbox{{\huge $0$}} \ \ \quad\quad \vdots \\
 & 1& & & a_1\\
 a_{n-1} & a_{n-2}& \dots & a_1 & \sum_{j=1}^{n-2} a_j a_{n-j-1} \\
\end{pmatrix}  (L^{\mathbb C}-\bar\lambda\cdot\Id)^n,
$$
\end{enumerate}
where
$$
\aligned
& a_1=  \lambda'_{z_n} z_1,  \\
& a_2=2  \lambda'_{z_n} z_2,  \\
& \dots                              \\
& a_{n-2}=(n-2) \lambda'_{z_n} z_{n-2}, \\
& a_{n-1}=1 + (n-1) \lambda'_{z_n} z_{n-1}.
\endaligned
$$

In the real coordinate system $x_1, y_1, x_2, y_2, \dots, x_n, y_n$  (where $z_k=x_k+iy_k$),  the operator $L$ and metric $g$ are defined by the $2n\times 2n$ real matrices which can be obtained from $L^{\mathbb C}$ and $g^{\mathbb C}$ by
the following standard rule:

--- each complex entry $a+ib$ of $L^{\mathbb C}$ is replaced by the $2\times 2$ block
$\begin{pmatrix} a & -b \\ b  & a \end{pmatrix}$,

--- each complex entry $a+ib$ of $g^{\mathbb C}$ is replaced by the $2\times 2$ block  $\begin{pmatrix}  a & -b \\ -b & -a \end{pmatrix}$.

Conversely, $g$ and $L$ defined by the above formulas are compatible for every holomorphic function $\lambda(z_n)$  (in the domain where $\lambda(z_n)\ne 0$ and $\det g\ne 0$, i.e. $1 + (n-1) \lambda'_{z_n} z_{n-1}\ne 0$).
\end{Th}

As we see, the case of a complex Jordan block is very similar to the real one. However, there is one very essential difference: the additional factor
$(L^{\mathbb C}-\bar\lambda\cdot\Id)^n$ in the formula for $g$.  Notice, by the way, that unlike $L^{\mathbb C}$ the components of $g^{\mathbb C}$ are not holomorphic functions on $M$ because of $\bar\lambda$ involved in this formula.

\begin{Rem}
If the differential of $\lambda$ does not vanish at the point $p$, then  just  as in the case of a real Jordan block,  we can
take $\lambda$ as the coordinate $z_n$ and obtain the following version of Theorem~\ref{onecomplexblock}.

Let $g$ and $L$ be compatible  and $L$ be conjugate to a complex Jordan block with complex conjugate eigenvalues $\lambda$ and $\bar\lambda$  {\rm (}$\mathrm{Im}\, \lambda\ne 0${\rm )}. If $d\lambda(p)\ne 0$,  then in a neighborhood of $p$ there exists a complex structure $J$ and a local complex coordinate system
$(z_1,\dots, z_n)$  on $M$ such that
\begin{enumerate}
\item
$L$ is a complex linear operator on $(T_PM, J)$ given in this coordinate system by the matrix:
$$
L^{\mathbb C} =
\begin{pmatrix}
z_n & 1            &       & &  z_1 \\
             & z_n & \ddots     & &  2z_2 \\
             &              & \ddots &  1 &  \vdots  \\
             &              &        & z_n & h(z_n)\!+\!(n\!-\!1)z_{n-1}         \\
             &              &        &              & z_n
\end{pmatrix}
$$
\item
$g$ is the real part of the complex bilinear form on $(T_PM, J)$ given in this coordinate system by the matrix:
$$
g^{\mathbb C}= -i
\begin{pmatrix}
 & & &   & h(z_n)\!+\!(n\!-\!1)z_{n-1} \\
 & & & 1 & (n-2) z_{n-2} \\
 & & \iddots & & \!\!\!\!\!\!\!\!\!\!\!\!\!\! \!\!\!\!\!\!\!\!\!\! \mbox{{\huge $0$}} \ \ \quad\quad \vdots \\
 & 1 &     & & z_1 \\
h(z_n)\!+\!(n\!-\!1)z_{n-1}& \dots & \dots &  z_1 & \sum
\end{pmatrix}(L^{\mathbb C}-\bar\lambda\cdot\Id)^n,
$$
where $\sum = \sum^{n-2}_{i=1} i (n-i+1) z_i z_{n-i-1}$
and $h(z_n)$ is a holomorphic function such that $h(z_n(p))\ne 0$.
\end{enumerate}

The passage from the complex coordinates $z_i$ to real coordinates $x_k, y_k$  ($z_k = x_k + i y_k$) follow the same rules as explained in Theorem \ref{onecomplexblock}.

\end{Rem}

\begin{Rem} We can easily rewrite this result in real coordinates. Namely, if the differential of the complex eigenvalue $\lambda$ does not vanish at the point $p\in M$, then in a neighborhood of this point 
there exists a local coordinate system $x_1, y_1, \dots, x_n , y_n$ such that the metric $g$ and operator $L$ are given as follows
$$
L= C^{-1}L_0 C,   \quad g=  C^\top g_0 \tilde L_0^n C
$$
where
$$
L_0=\begin{pmatrix}
Z_n & \Id_2 &   &   &  \\
        & Z_n &  \ddots & & \\
        &        & \ddots &\ddots &\\
        &   & & Z_n & \Id_2 \\
        & & & & Z_n
\end{pmatrix}, 
\qquad
\tilde L_0=\begin{pmatrix}
2iY_n & \Id_2 &   &   &  \\
        & 2iY_n &  \ddots & & \\
        &        & \ddots &\ddots &\\
        &   & & 2iY_n & \Id_2 \\
        & & & & 2iY_n
\end{pmatrix},
$$
$$
g_0=\begin{pmatrix}
 & & & & \tilde{ \bf 1}_2 \\
 & & &  \tilde{ \bf 1}_2 &   \\
 & &\iddots & & \\
 & \tilde{ \bf 1}_2 & & &    \\
 \tilde{ \bf 1}_2 & & & &
\end{pmatrix}, \quad 
C=\begin{pmatrix}
\Id_2 & & & & {\bf 0}_2\\
 & \Id_2 & & &  Z_1\\
&  & \ddots & & \vdots \\
& & & \Id_2 & (n-2) Z_{n-2}\\
& & & & H+(n-1)Z_{n-1}
\end{pmatrix}. 
$$
Each of the indicated entries  represents a $2\times 2$-matrix of the following form: \\
$
\Id_2 = \begin{pmatrix} 1 & 0 \\ 0 & 1 \end{pmatrix}$, $\tilde{ \bf 1}_2= \begin{pmatrix} 0 &  1\\  1 & 0 \end{pmatrix}$,  ${\bf 0}_2 = \begin{pmatrix}  0&  0\\   0& 0 \end{pmatrix}$, $Z_i= \begin{pmatrix}  x_i & -y_i \\  y_i & x_i \end{pmatrix}$,  $2iY_n =  \begin{pmatrix}  0 & -2y_n \\  2y_n & 0 \end{pmatrix}$ and $H= \begin{pmatrix} u & -v  \\ v  & u \end{pmatrix}$ with $u=u(x_n,y_n)$, $v=v(x_n,y_n)$ being functions satisfying the Cauchy--Riemann conditions  (i.e.,  $h=u+i v$ is a holomorphic function of $z_n=x_n+iy_n$).
\end{Rem}




 \subsection{ Perspectives and first global results}


It is hard to overestimate the  role of the Levi-Civita theorem in the local and global theory of  geodesically equivalent  Riemannian metrics. Almost all local results are based on it or can easily be proved using it.  Though Levi-Civita theorem is local,  most global (i.e., when the manifold is compact) results on geodesically equivalent Riemannian metrics  also use Levi-Civita theorem as an important tool. Roughly speaking, using the Levi-Civita description one can reduce  any problem that can be stated using geometric PDEs (for example, any problem involving the curvature) to solving or analysing  a system of ODEs.

We expect that our result will play the same role in the pseudo-Riemannian case. We suggest using it to prove the natural generalization of the projective  Lichnerowicz-Obata conjecture  and the Sophus Lie  problem for the pseudo-Riemannian case, see \cite[\S2.2]{splitting} for the description of the problems. Note that the Lichnerowicz-Obata conjecture was solved, in the Riemannian case, in  \cite{sol} under additional assumptions and in \cite{archive} in full generality, and the solution essentially used the Levi-Civita theorem. The Sophus Lie  problem was solved in the Riemannian case for dimensions $n>2$
 in  \cite{sol}, the solution again used the Levi-Civita theorem
 and in  the 2-dimensional case for all signatures  in \cite{alone} where the solution essentially used the description of two dimensional projectively equivalent Riemannian and pseudo-Riemannian metrics obtained for example in   \cite{appendix}.

We also hope that our description will be helpful in understanding of the global structure of the manfolds carrying geodesically equivalent pseudo-Riemannian metrics. One of the ultimate goals could be to understand the ``possible topology'' of such manifolds. Though our main  theorem is local, it can be effectively used (as was the case with the Levi-Civita theorem)  in the global setting too.
In particular, we prove the following two results.

\begin{Th}
\label{appl1}
Let $M$ be a closed  connected  manifold. Suppose  $g$ and $\bar g$ are geodesically equivalent metrics on it and $L$ given by \eqref{L} has a complex eigenvalue $\lambda\in \mathbb C\setminus \R$ at least at one point.
Then at every point of $M$ this number $\lambda$ is an eigenvalue of $L$. Moreover, the multiplicity of the eigenvalue $\lambda$ is the same at every point of the manifold. 
\end{Th}

\begin{Cor} \label{appl2}
Let $M^3$ be a closed connected 3-dimensional manifold. Suppose  $g$ and $\bar g$ are geodesically equivalent metrics on it and at least at one point, $L$ given by \eqref{L} has a non real eigenvalue.
Then $M^3$ can be finitely covered by the 3-torus.
\end{Cor}



   \section{Proof of Theorem \ref{multtwoth}: case of geometric multiplicity {\mathversion{bold}$\ge 2$}}
\label{affine}


We assume that $(M, g)$ is connected and  that (a selfadjoint $(1,1)$-tensor)  $L$ is compatible with $g$. Our first goal is to prove  

\begin{proposition}\label{multtwo}
  Assume that  in a neighborhood  $U\subseteq M$  there exists a  continuous  function $\lambda:U\to \mathbb{R}$ or $\lambda:U\to \mathbb{C}$  such that for every $x\in U$  the number $\lambda(x)$  is an eigenvalue of $L$  at $x$ of geometric multiplicity at least two.     Then,   the function $\lambda$ is  constant; moreover, for every point 
    $x\in M$  the number $\lambda$ is an eigenvalue  of $L$ at $x$ of geometric multiplicity  at least two.
\end{proposition}

{\it Proof.}   Our proof will use the following theorem due to \cite{benenti,MT,topalov}. For any $(1,1)$-tensor
   $A$ on $M$, let us denote by $\operatorname{co} (A)^\top$ the $(1,1)$-tensor whose matrix in a local coordinate system 
 is     the comatrix of  (= adjoint  matrix to) $A$ transposed. It is indeed  a well-defined tensor field: smoothness follows from the fact that the components of $\operatorname{co}(A)^\top$ are algebraic expressions in the entries of $A$. The change-of-basis transformation law holds for $\operatorname{co} (A)^\top$,  if $A$ is nondegenerate, since in this case $\operatorname{co}(A)^\top = \det(A) A^{-1}$. As nondegenerate matrices are dense in the set of all quadratic matrices, the transformation law holds for any $A$.

   \begin{Th} Let $L$ be compatible with $g$.  Then for any $t\in \mathbb{R}$, the function
      \begin{equation} \label{integral} I_t:TM\to \mathbb{R},  \quad  I_t(\xi)=  g(\operatorname{co}(L- t\cdot \Id)^\top \xi, \xi)
  \end{equation} is an integral of the geodesic flow of the metric $g$.
   \end{Th}

   Recall that a function $I$ is an {\it integral}, if for every geodesic $\gamma$  parameterized by a natural  parameter~$s$ (such that $\nabla_{\dot \gamma}\dot \gamma = 0$),     the function $s\mapsto I(\dot \gamma(s))$ is constant.

We first consider the case when the function $\lambda$ is real. 
As the algebraic multiplicity of the eigenvalue $\lambda$ is upper semi-continuous, replacing perhaps $U$ by an open subset of it, we may assume that  the algebraic multiplicity  is constant on $U$.  
Then,  the implicit function theorem easily implies that $\lambda$ is a smooth function.

   First we prove that $\lambda$ is constant on $U$.
    By contradiction, assume that there exists   $p\in U$ where the differential of $\lambda$ is  not zero.  Then, in a small neighborhood of $p$, the set
$$
M_{\lambda(p)}:= \{ q\in M \mid \lambda(q)= \lambda(p)\}$$
is a smooth submanifold of $M$ of codimension $1$.
At every  point of $M_{\lambda(p)}$,  the matrix of  the tensor $(L- \lambda(p)\cdot \Id)$ has rank at most $n-2$, so
$\operatorname{co}(L- \lambda(p)\cdot \Id)^\top = { \textrm{\bf\large 0}} $. Consequently, for every point $q\in M_{\lambda(p)}$
 and for every $\xi\in T_qM$ we have $I_{\lambda(p)}(\xi)= 0$.
 Now, take a point $x\in U$, 
 $x \not\in M_{\lambda(p)}$, and consider all geodesics $\gamma_{q,x}$ connecting the points  $q\in M_{\lambda(p)}$ with $x$. We assume  that the parameter $s$ on the geodesic is natural, $\gamma_{q,x}(0)=q\in  M_{\lambda(p)}$   and   $\gamma_{q,x}(1)= x$.  Since $ M_{\lambda(p)} $ has codimension one, for all 
  $x$ that are  sufficiently close to $p$,  the set of vectors that are proportional to the  velocity  vectors $\dot\gamma_{q,x}(1)$ of such geodesics,  contains an open nonempty subset of $T_xM$.  Then, 
 $ \operatorname{co}(L- \lambda(p)\cdot \Id)^\top $ is zero at the point $x$.  It follows immediately  that  $\lambda(p)$ is an eigenvalue of $L$ at the point $x$.
 Then, $\lambda$ is constant in a neighborhood of $p$  which contradicts our assumption that $d\lambda_{|p}\ne 0$.  The contradiction shows that $\lambda$ is constant on $U$.

 Let us now consider the case when   $L$  has two complex conjugate eigenvalues $\lambda,\bar \lambda: U\subseteq M\to \mathbb{C}$, $\mathrm{Im}\, \lambda \ne 0$.  We again assume without loss of generality that
  the algebraic multiplicity of the eigenvalue $\lambda(x)$ 
  is the same at all
points $x \in  U$, which in particular implies  that  $\lambda$  is a smooth function.  We first note that, for every $(1,1)$-tensor $A$,  the $(1,1)$-tensor
   $\operatorname{co}(A- t \cdot \Id)^\top$ is a polynomial in $t$ of degree $n-1$ whose coefficients are $(1,1)$-tensors.
    Then,  for every complex number $\tau$, 
    the real and  imaginary parts of the complex-valued function
    $$I_\tau :TM\to \mathbb{C}, \    I_{\tau}(\xi):=   g(\operatorname{co}(L- \tau \cdot \Id)^\top \xi, \xi)$$
   are also integrals.  Since $\mathrm{rank}\,( L(q)-\lambda(q)\cdot\Id )\le n-2$, for every $q$ such that $\lambda(q)= \tau$ we have that
      $I_\tau(\xi)= 0$ for every $\xi\in T_qM$.

 Suppose $\lambda$ is not constant on $U$. Then for a certain point $p$  of  $U$ its differential is not zero.   Suppose first that the differential of the real part of $\lambda$ is proportional to the differential of the imaginary part  at all points of a certain neighborhood of $p$. Then, in a sufficiently small neighborhood $U'(p)\subseteq  U$ of
$p$ the set 
$M_{\lambda(p)}:= \{ q\in M \mid \lambda(q)= \lambda(p)\}$ 
is a submanifold of dimension $n-1$, as it was in the case of  a real eigenvalue $\lambda$.
 Then, repeating  the same arguments as above we conclude that $\lambda(x)= \lambda(p)$ for all $x$ from a small neighborhood of $p$, which gives us a contradiction with the assumption that the differential of $\lambda$ does not vanish at $p$. The contradiction shows that $\lambda$ is a constant   provided  the differential of the real part of $\lambda$ is proportional to the differential of the imaginary part in some $U'\subseteq U$.

Let us now suppose that the differential of the real part of $\lambda$ at the point $p$  is  not proportional to the differential of the imaginary part. Then, the set $
M_{\lambda(p)}:= \{ q\in M \mid \lambda(q)= \lambda(p)\}$ is (in a sufficiently small neighborhood $U'(p)\subseteq U$) a submanifold of dimension $n-2$. We again take an arbitrary point $x$ that is sufficiently close to $p$ and consider all geodesics $\gamma_{q,x}$ connecting the points $q\in M_{\lambda(p)}$ with  $x$ assuming as above that
  $\gamma_{q,x}(0)\in  M_{\lambda(p)}$ and   $\gamma_{q,x}(1)= x$. The set of the  tangent
  vectors  at $x$ that are proportional to the velocity vectors  to such geodesics  at the point $x$
  contains a submanifold of codimension  1 of $T_xM$ implying that the real part of $ I_{\lambda(p)}(\xi)$ is proportional to the imaginary part of $ I_{\lambda(p)}(\xi)$  for all $\xi \in T_xM$ (the coefficient of the proportionality is a constant on each $T_xM$ but may a priori depend on $x$). Now, since both functions, the  real and the imaginary parts of $ I_{\lambda(p)}$,  are integrals, the coefficient of proportionality of these functions is an integral too implying it is constant. Then,  for
  a  certain complex  constant $a+ i b\ne 0$, for every $x\in U'(p)$ and every $\xi\in T_xM$ we have   
  $(a + i b)  I_{\lambda(p)}(\xi) = (a - i b)I_{\bar \lambda(p)}(\xi)$ so that
 \begin{equation} 
 \label{mat1}  
 (a+ib) \operatorname{co}(L- \lambda(p)\cdot  \Id  )^\top=  (a-ib) \operatorname{co}(L- \bar \lambda(p) \cdot \Id  )^\top .
 \end{equation} 
   For points $x$ such that $\lambda(x)\ne \lambda(p)$ the matrix $L- \lambda(p) \cdot \Id  $ is nondegenerate   and  \eqref{mat1} implies  that $(L- \lambda(p) \cdot \Id)^{-1}  $ is proportional to $(L- \bar\lambda(p) \cdot \Id )^{-1}$.   Hence, $L- \lambda(p) \cdot \Id$ and $L- \bar\lambda(p) \cdot \Id$ are proportional too, which  contradicts the assumption that $\lambda$ is not real. The contradiction shows that at all points of the neighborhood $U$ (such that $\lambda(x)$ is a eigenvalue of $L$ of geometric multiplicity at least two at every point $x$) the function $\lambda $ is a constant.

 Let us now show that this (real or not) constant $\lambda$  is an eigenvalue of $L$ of geometric  multiplicity at least two at every point of the  whole $M$.    We
 first  consider a point $p\in M\setminus U$  that can be connected 
  with a point of  $U$ by a geodesic $\gamma$,  where $U$ is a neighborhood such that at  each its point $L$ has (constant) eigenvalue $\lambda$ of multiplicity at least 2; we think that $\gamma(0)  =p$ and $\gamma(1)\in U$. We consider a small open 
   neighborhood $V\subseteq   T_pM$ of $\xi= \dot\gamma(0)$. If $V$ is sufficiently small, 
   for every $\eta \in V$  the point  $ \gamma_{p, \eta} (1)$ of the  geodesic $\gamma_{p, \eta}$ such that $\gamma(0)= p$ and $\dot \gamma(0)= \eta$   lies in $U$.  
   Since at each point of $U$ the constant  $\lambda$ is an eigenvalue of $L$ of multiplicity at least 2,  $ I_{\lambda} (\dot \gamma_{p, \eta} (1))=0$
   implying $ I_\lambda(\dot \gamma_{p, \eta} (0))=0$. 
   Hence,   $I_\lambda\equiv 0$ on an open nonempty  subset of $T_pM$  implying   $I_\lambda\equiv 0$  on  the whole $T_pM$  so $\lambda$ is an eigenvalue of $L$ at $p$ of multiplicity at least two. Now, if $p$ can be connected  by a geodesic with a point of $U$, then any point from a sufficiently small neighborhood of $p$ can also be connected by a  geodesic with a point of $U$ so $\lambda$ is an eigenvalue of $L$  of multiplicity at least two at every point of a small neighborhood of $p$. 
  To come to the same conclusion  on the whole $M$, it suffices to notice that every point $x\in M$ can be joined with $U$ by a piecewise smooth  curve such that  each smooth segment of it  is a geodesic. 
   Proposition \ref{multtwo} is proved.  \qed
\begin{Cor}
In the hypotheses of Proposition \ref{multtwo},  assume in addition that $\lambda\in\R$ is the unique eigenvalue of $L$  (or  $\lambda, \bar\lambda\in \mathbb C\setminus \R$ are the unique eigenvalues of $L$), then $L$ is covariantly constant, i.e.,
$\nabla L=0$ or,  equivalently, $g$ and $\bar g$ are affinely equivalent.  
\end{Cor}

The proof is obvious:  since $\lambda$ is constant, so is $\tr L$.  Hence, the right hand side of \eqref{main} vanishes and we get $\nabla L=0$. \qed

Thus, if the eigenspace of $L$ has dimension $\ge 2$,  then our problem is reduced to the classification of pairs of affinely equivalent pseudo-Riemannian metrics, which  
has been recently obtained by Boubel in  \cite{boubel}.



\section{Proof of Theorem \ref{onerealblock}: case of a real Jordan block}
\label{proof2}


\subsection{Canonical frames and uniqueness lemma}
\label{frame}


Let $L$ be $g$-selfadjoint operator on a real vector space $V$.  It is a natural question to ask to which canonical form we can  simultaneously reduce (the matrices of) $L$ and $g$ by an appropriate change of a basis. The answer is well known (see, for example, \cite{lancaster}) and is given by the following

\begin{proposition}
\label{prop2}
There exists a {\it canonical} basis $e_1, \dots, e_n\in V$ in which $L$ and $g$ can be simultaneously reduced to the following block diagonal canonical forms:
$$
L_{\mathrm{can}}=
\begin{pmatrix}
L_1 &         & &  \\
        & L_2 & &  \\
        &         & \ddots &  \\
        &         &             & L_s
\end{pmatrix} ,
\quad
g_{\mathrm{can}}=
\begin{pmatrix}
g_1 &         & &  \\
        & g_2 & &  \\
        &         & \ddots &  \\
        &         &             & g_s
\end{pmatrix} ,
$$
where
\begin{equation}
\label{realJordan}
L_j =
\begin{pmatrix}
\  \lambda & 1 & &  \\
                        &\lambda & \ddots &  \\
                        &                        & \ddots & 1  \\
                        & & &  \lambda \
\end{pmatrix}
\end{equation}
in the case of a real eigenvalue $\lambda\in \R$ (real Jordan block), or
\begin{equation}
\label{complexJordan}
L_j=
\begin{pmatrix}
\begin{matrix} a & \!\!\! -b \\ b & a \end{matrix}   &  \begin{matrix} 1 & 0 \\ 0 & 1 \end{matrix} & & \\
  & \begin{matrix} a & \!\!\! -b \\ b & a \end{matrix} & \ddots  & \\
  &                  &    \ddots  &  \begin{matrix} 1 & 0 \\ 0 & 1 \end{matrix} \\
  &    &    & \begin{matrix} a & \!\!\! -b \\ b & a \end{matrix}
  \end{pmatrix}
\end{equation}
  in the case of complex conjugate eigenvalues $\lambda_{1,2}=a\pm ib$, $b\ne 0$ (complex Jordan block), and
  \begin{equation}
\label{gcanon}
g_j =
\varepsilon \ \begin{pmatrix}
 &  & &   & 1 \\
 &  & & 1 &   \\
 &  & \iddots & & \\
 & 1 &      &  & \\
1 &  &      &  &
\end{pmatrix},
\end{equation}
where $\varepsilon=\pm1$  in the case $\lambda\in\R$ and $\varepsilon = 1$ for $\lambda_{1,2}= a\pm ib\in\mathbb C\setminus \R$.
It is assumed that for each $j$ the blocks $g_j$ and $L_j$ are  of the same size and  that the corresponding eigenvalues depend on $j$.
\end{proposition}

\begin{Rem}\label{addremark}
Notice that the canonical forms $g_{\mathrm{can}}$ and $L_{\mathrm{can}}$ can be chosen in many different ways.  For example, in the complex case we can replace $g_{\mathrm{can}}$ by $g_{\mathrm{can}} P(L_{\mathrm{can}})$ where $P(t)$ is an arbitrary polynomial such that $P(L_{\mathrm{can}})$ is invertible.  Indeed, as $L_{\mathrm{can}}$ is selfadjoint w.r.t. $g_{\mathrm{can}}$,  the pair $(g_{\mathrm{can}} P(L_{\mathrm{can}}),  L_{\mathrm{can}})$ is also a pair  consisting of a nondegenerate symmetric bilinear form and a selfadjoint operator w.r.t. it.  Moreover, as $L_{\mathrm{can}}$ is ``complex'', it follows from Proposition \ref{prop2} that 
$(g_{\mathrm{can}}, L_{\mathrm{can}})$  and $(g_{\mathrm{can}} P(L_{\mathrm{can}}),  L_{\mathrm{can}})$ are conjugate to each other,   since in the ``complex'' case the canonical form for $g$ is uniquely defined by $L$. For our purposes, by a canonical form of $L$ and $g$ it is convenient to understand any forms where the entries of $L_{\mathrm{can}}$ and $g_{\mathrm{can}}$ depend on the eigenvalues of  $L$ only.
\end{Rem}

Now let $g$ be a pseudo-Riemannian metric on a smooth manifold $M$ and $L$ be a $g$-selfadjoint $(1,1)$-tensor field.  Assume that $L$ is regular at each point of a small neighborhood $U(p)$ of a point $p\in M$. Recall that the {\it regularity} of $L$ means that each eigenvalue of $L$  is of geometric multiplicity one or, equivalently, the Jordan normal form contains exactly one Jordan block for each eigenvalue.  This condition implies that the eigenvalues of $L$ are smooth functions on $U(p)$ and
the Jordan type of $L$ does not change in $U(p)$, in particular, $p$ is a regular point (see page \pageref{vocal}).
In such a situation we can choose smooth linearly independent vector fields $e_1, \dots, e_n\in T_xM $ in which
$L$ and $g$ both take canonical forms. For a regular $L$, these vectors are uniquely defined  (up to a discrete group) and we will say that $e_1, \dots, e_n$ is a {\it canonical moving frame} for $L$ and $g$.

In general, the vector fields $e_1, \dots, e_n$ do not commute. To reconstruct a canonical coordinate system on $U(p)$ we need to analyse the commutation relations between them. It turns out that these relations can be obtained from the compatibility equation \eqref{main}.

\begin{Lemma}
\label{uniqueness}
Let
$e_1, \dots, e_n$ be a canonical moving frame for $L$ and $g$ in a neighborhood of a regular point $p\in M$.
If $L$ and $g$ are compatible and $L$ is regular, then the covariant derivatives $\nabla_{e_i} e_j$ and hence
the commutators $[e_i, e_j]$ can be uniquely expressed as certain linear combinations of
$e_l$ with the coefficients being functions of $\lambda_r$ and their derivatives $e_s(\lambda_r)$ along $e_s$, where $\lambda_r$ are eigenvalues of $L$.
\end{Lemma}

{\it Proof.}
For the frame $e_1, \dots, e_n$ we introduce $B_u$ to be $(1,1)$-tensor field defined by
\begin{equation}
\label{Bu}
B_u v = \nabla_u v,
\end{equation}
where $u$ and $v$ are vector fields with constant coordinates w.r.t. the frame.

Clearly, $B_u$ defines the Levi-Civita connection in the frame $e_1, \dots, e_n$ and our goal is to reconstruct it from the compatibility equation \eqref{main}. The covariant derivative of $L$ in terms of $B_u$ can be written as
$$
\nabla_u L = \mathcal D_u(L) + [B_u, L]
$$
where $\mathcal D_u L$ denotes the operator obtained by differentiating $L$ componentwise along $u=u^k e_k$,  i.e., for $L=L^i_j e_i\otimes e^j$, we have $\mathcal D_u L= u^k e_k(L^i_j) e_i\otimes e^j$.

To find $B_u$, it is convenient to rewrite the compatibility equation in the form
\begin{equation}
\label{a0}
[B_u, L] = \frac{1}{2} (u\otimes d\tr L + (u\otimes d\tr L )^* ) - \mathcal D_u L.
\end{equation}

In addition to that we have
$$
\mathcal D_u \bigl(g(e_i, e_j)\bigr) = \nabla_u \bigl( g(e_i, e_j) \bigr)=
g(B_u e_i, e_j) + g(e_i, B_u e_j)
$$
or, equivalently
$$
B_u + B_u^* = g^{-1} \mathcal D_u g
$$

Thus, $B_u$ satisfies two equations of the form
\begin{equation}
\begin{array}{ll}
[B_u, L] = C \\
B_u + B_u^* = D
\label{a6}
\end{array}
\end{equation}
where $C$ and $D$ are certain operators  (whose components w.r.t. the moving frame are functions of the eigenvalues $\lambda_r$ and their derivatives $e_s(\lambda_r)$).

The uniqueness of the solution  (if it exists!) is a purely algebraic fact.
Indeed, consider the corresponding homogeneous system
\begin{equation}
\begin{array}{ll}
[B_u, L] = 0 \\
B_u + B_u^* = 0
\label{a7}
\end{array}
\end{equation}

The first equation means that $B_u$ commutes with $L$, i.e., belongs to the centralizer of $L$.  Since $L$ is regular,  its centralizer is generated by powers of $L$,
i.e., $\Id, L, L^2, \dots, L^{n-1}$.  It follows from this that $B_u$ is $g$-selfadjoint. But then the second equation can be rewritten
simply as $2B_u = 0$.  Thus, the homogeneous system has only the trivial solution which proves the statement.

The commutators $[e_i, e_j]$ can now be uniquely reconstructed
by means of the standard formula: $[e_i, e_j]=\nabla_{e_i} e_j - \nabla_{e_j} e_i = B_{e_i} e_j - B_{e_j} e_i$.
\qed

In the  next section we show how the commutation relations between the elements of the canonical moving frame can be found in practice.

Notice that Lemma \ref{uniqueness} does not say that \eqref{a6} is always consistent for every $C$ and $D$. In fact, these matrices have to satisfy  some additional relations (for example, $\tr CL^k=0$). These equations, in particular, imply the vanishing of the Nijenhuis torsion of $L$ and, therefore, the fact that $e_s(\lambda_r)=0$ for those $e_s$ which ``do not belong'' to the $\lambda_r$-block. Another condition of this kind is discussed below in Lemma \ref{en}.


\subsection{Canonical frame and canonical coordinate system for a real Jordan block}
\label{real}


Let $L$ be conjugate to a Jordan block with a real eigenvalue. Then we can choose a moving frame $e_1,\dots, e_n$ in which $L$ and $g$ take the following canonical forms:
\begin{equation}
\label{canonicalJordan}
L_{\mathrm{can}} =
\begin{pmatrix}
\lambda(x) & 1 & &  \\
                        &\lambda(x) & \ddots &  \\
                        &                        & \ddots & 1  \\
                        & & &  \lambda(x)
\end{pmatrix}, \qquad
g_{\mathrm{can}} = \pm \begin{pmatrix}
 &  & & 1 \\
                        & & 1 &  \\
                        &      \iddots                  &  &   \\
                     1   & & &
\end{pmatrix}
\end{equation}

Here we apply the ideas from \S\,\ref{frame} to describe the commutation relations between $e_1,\dots, e_n$ and then to solve them in order to construct a canonical coordinate system for the pair $g$ and $L$.

As usual, it is convenient to decompose $L$ canonically into the semisimple and nilpotent parts:
$$
L=\lambda(x)\cdot \Id + N.
$$

Obviously, $N$ is selfadjoint with respect to $g$. The compatibility equation can naturally be rewritten in terms of $N$:
$$
\nabla_u(L)= u(\lambda)
\cdot \Id + \nabla_u N = \frac{n}{2} \bigl( u\otimes d\lambda + (u\otimes d\lambda)^* \bigr)
$$
or
\begin{equation}
\label{nil1}
[B_u, N] = \frac{n}{2} \bigl( u\otimes d\lambda + (u\otimes d\lambda)^*\bigr ) - u(\lambda) \cdot \Id,
\end{equation}
where $B_u$, as before, is defined by \eqref{Bu} and we use the fact that the components of $N$ in the frame are all constants so that
$\nabla_u N = [B_u, N]$.
This equation implies

\begin{Lemma}
\label{en}
We have $e_i(\lambda)=0$, for $i=1,\dots, n-1$.
\end{Lemma}

{\it Proof.}
This property is well known for $L$ with zero Nijenhuis torsion  (for $L$ this condition is fulfilled, see e.g. \cite[Theorem 1]{benenti}). However, we can easily derive this fact from \eqref{nil1}.  Indeed, multiply both sides of this equation by $N$ and take the trace. For the left hand side we get:
$$
\tr \bigl( N\cdot [B_u,N] \bigr) = 0.
$$
 For the right hand side:
$$
\tr  \left(N \cdot \Bigl(\frac{n}{2} \bigl(u\otimes d\lambda + (u\otimes d\lambda )^* \bigr) - u(\lambda)\cdot \Id \Bigr)\right)=
n\cdot N\! u \,(\lambda) - u(\lambda)\cdot \tr N= n \cdot N\!u\,(\lambda).
$$
Hence, for any vector $v=N\!u \in \Image{N}$, we have $v(\lambda)=0$. It remains to notice that $\Image N = \mathrm{span} (e_1,\dots, e_{n-1})$. \qed

\begin{Rem}  Since this property plays an important role in many other problems appearing in differential geometry, we give the proof of Lemma \ref{en} under the only assumption that $N_L=0$. It was shown in \cite[Lemma 1]{benenti},  that  this assumption implies the following identity for $L$:
$$
d\,\trace L \, (u)=(d\,\ln \det L)(Lu)  \ \ \mbox{or,  equivalently,} \ \   u\, (\trace L) = Lu\, (\ln \det L)
$$
where $u$ denotes an arbitrary vector field.
 In our case $\trace L= n\lambda$ and $\det L=\lambda
^n$.  Therefore, 
$$
n u\, (\lambda) = n \lambda^{-1} Lu \, (\lambda) 
$$
and, consequently, 
$$ 
Lu (\lambda)   - \lambda  u(\lambda) = \bigl((L-\lambda\cdot \Id)u\bigr)  (\lambda)= N\!u \, (\lambda)=0. 
$$ 
In other words,
$v(\lambda)=0$ for any $v$ of the form $v=N\!u\in  \Image{N}$, as required.
\end{Rem}

Thus, in our basis $d\lambda=(0,\dots, 0, e_n(\lambda))$.
This allows us to get the following explicit form for the right hand side of \eqref{nil1}:
\begin{equation}
\label{forBu}
\nabla_u N= [B_u, N]= e_n(\lambda)
\begin{pmatrix}
\frac{n-2}{2}u_n & \frac{n}{2} u_{n-1} & \frac{n}{2} u_{n-2} & \dots  &
\frac{n}{2} u_{2} & n u_1 \\
                 & - u_n               &       & &                  &
\frac{n}{2} u_{2}\\
                 &                     & -u_n  &   \quad  \mbox{{\huge $0$}} &                  &
\frac{n}{2} u_{3}\\
                 &                     &       & \ddots &           &  \vdots \\
                 &                     &       &        &  -u_n     &
 \frac{n}{2} u_{n-1} \\
                 &                     &       &        &           &
\frac{n-2}{2}u_n
\end{pmatrix}
\end{equation}

According to Lemma~\ref{uniqueness}, page~\pageref{uniqueness}, the solution of this equation is unique.
We just give the final answer (the reader can check this result by
substituting \eqref{eq2} into \eqref{forBu}).
\begin{equation}
B_u =  e_n(\lambda)
\begin{pmatrix}
\frac{n}{2}u_{n-1} & \frac{n}{2}u_{n-2} & \dots & \frac{n}{2}u_{2}
&\frac{n}{2}u_{1} &  0 \\
(1-\frac{n}{2})u_n  &                    &
& & & -\frac{n}{2}u_{1}  \\
                   & (2-\frac{n}{2})u_n  &     &  & & -\frac{n}{2}u_{2}  \\
                   &                    & \ddots& &   \!\!\!\!\!\!\!\!\!\!\!  \mbox{{\huge $0$}}   &      \vdots             \\
                   &                     &     & (\frac{n}{2}-2)u_n & &
-\frac{n}{2}u_{n-2} \\
                   &           &         &       & (\frac{n}{2}-1)u_n &
-\frac{n}{2}u_{n-1}
\end{pmatrix}
\label{eq2}
\end{equation}

The next step is to find pairwise commutators $[e_i,e_j]$.

\begin{Lemma}
\label{firstcommutators}
The vector fields $e_1,\dots,e_{n-1}$ commute.
\end{Lemma}

{\it Proof.} Let $u=u_1e_1+\dots u_ne_n$.
It follows from \eqref{eq2} that
$$
\nabla_u e_j=e_n(\lambda) \left(\frac{n}{2} u_{n-j} e_1 + (j-\frac{n}{2}) u_n e_{j+1}\right), \quad
j<n.
$$
Hence, for $i<n$, $\nabla_{e_i} e_j=e_n(\lambda)\frac{n}{2} e_1$ if and only if $i+j=n$,
otherwise $\nabla_{e_i} e_j=0$. In any case $$ [e_i,e_j]= \nabla_{e_i} e_j
-\nabla_{e_j} e_i =0 $$ for $i,j<n$.
\qed

It remains to find the commutators $[e_i,e_n]$.

\begin{Lemma}
\label{secondcommutators}
For $i=1,\dots, n-1$, we have
$$
[e_i,e_n]=- i  e_n(\lambda) \cdot e_{i+1}.
$$
\end{Lemma}

{\it Proof.} From \eqref{eq2} we have
$$
\nabla_{e_i} e_n = - \frac{n}{2}  e_n(\lambda)\cdot e_{i+1} \quad
\mbox{and}
\quad
\nabla_{e_n} e_i = -(i- \frac{n}{2}) e_n(\lambda)\cdot e_{i+1}.
$$
Thus,
$$
[e_i,e_n]=- e_n(\lambda)\left( \frac{n}{2} e_{i+1}+ -(i- \frac{n}{2})
e_{i+1}\right)=- i  e_n(\lambda)\cdot e_{i+1},
$$
as stated.  \qed

Our goal now is to find a coordinate system with respect to which $N$
and $g$ have the simplest form. Since the vector fields $e_1,\dots, e_{n-1}$
commute we can choose a coordinate system $x_1,\dots,x_n$ in such a way
that $e_1=\partial_{x_1},\dots,e_{n-1}=\partial_{x_{n-1}}$. To make our choice unambiguous, we assume that our initial point $p\in M$ has all coordinates zero and, in addition,
\begin{equation}
\label{initial}
e_n=\partial_{x_n} \quad \mbox{on the $x_n$-axes,}
\end{equation}
i.e. on the curve $x_1=x_2=\dots=x_{n-1}=0$. Notice that
the foliation generated by $\Image N$ is given by $x_n=const$ and the
eigenvalue $\lambda$ depends on $x_n$ only.

To rewrite $L$ and $g$ in this coordinate
system we just need to find the transition matrix between $e_1,\dots, e_n$ and $\partial_{x_1},\dots ,\partial_{x_n}$. Since $\partial_{x_i}=e_i$, $i=1,\dots, n-1$, it remains to determine the coefficients (yet unknown) of the linear combination
$$
\partial_{x_n}= a_0e_1 + \dots + a_{n-1} e_n
$$

First we use the fact that $\lambda$ does not depend
on $x_1,\dots,x_{n-1}$. Therefore
$$
\lambda'_{x_n}=\partial_{x_n}(\lambda)=
(a_0e_1 + \dots + a_{n-1} e_n)(\lambda) =a_{n-1} e_n(\lambda).
$$

Since $\partial_{x_n}$ must
commute with each $e_i=\partial_{x_i}$ ($i<n$),
we obtain a system of differential equations on $a_j$:
$$
0=[e_i, a_0e_1 + \dots +a_{n-1} e_n]= \sum_{l=1}^n \frac{\partial a_{l-1}}{\partial
x_i}\cdot e_l - a_{n-1} \, i e_n(\lambda)\cdot e_{i+1} = \sum_{l=1}^n \frac{\partial a_{l-1}}{\partial
x_i} \cdot e_l -  i \lambda'_{x_n} \cdot e_{i+1} ,
$$
or, equivalently,
$$
\frac{\partial a_{l-1}}{\partial
x_i}=0, \quad \hbox{if \ \ $l\ne i+1$ and \ \ }
\frac{\partial a_{i}}{\partial
x_i}= i \lambda'_{x_n}, \quad i=1,\dots, n-1.
$$

In other words, $a_0=a_0(x_n)$, whereas  $a_{i}$ depends on $x_{i}$ and
$x_n$ and satisfies the equation
$$
\frac{\partial a_{i}}{\partial
x_{i}}= i\, \lambda'_{x_n},\quad   i = 1,\dots, n-1,
$$
which can be easily solved. Its general solution is
$$
a_{i}(x_{i},x_n)= i\, \lambda'_{x_n} x_{i} +
f_i(x_n),
$$
where $f_i(x_n)$ is an arbitrary function. But we have a kind of {\it initial condition} \eqref{initial} that requires
$$
a_{i}(0,\dots, 0, x_n)=0 \quad \mbox{for } i\ne n-1, \quad \mbox{and}\quad a_{n-1}(0,\dots, 0, x_n)=1.
$$
It follows immediately from this that
\begin{equation}
\label{ai}
\begin{aligned}
a_0 &=0 \\
a_{i}& = i\, \lambda'_{x_n} x_{i}, \quad i=1,\dots, n-2,
\end{aligned}
\end{equation}
and
\begin{equation}
\label{an}
a_{n-1}=(n-1)  \lambda'_{x_n} x_{n-1} +
1.
\end{equation}

Thus, the transition matrix $C$ has been found:
$$
(\partial_{x_1}, \dots, \partial_{x_n}) = (e_1,\dots, e_n) \cdot C, \quad \mbox{with} \ \
C=\begin{pmatrix} 1 & & & a_0 \\ & 1 & & a_1 \\ & & \ddots & \vdots \\ & & & a_{n-1} \end{pmatrix}
$$
and $a_0, \dots, a_{n-1}$ defined by \eqref{ai} and \eqref{an}.

Now to obtain the form of $g$ and $L$ in the local coordinates $x_1,\dots, x_n$, we only need to apply the standard rule
$$
L_{\mathrm{can}}  \longrightarrow L= C^{-1} L_{\mathrm{can}} C, \quad  g_{\mathrm{can}} \longrightarrow  g=C^\top g_{\mathrm{can}}C.
$$
where $L_{\mathrm{can}}$ and $g_{\mathrm{can}}$ are defined by \eqref{canonicalJordan}.
A straightforward computation of $L$ and $g$ gives the statement of Theorem \ref{onerealblock}, page \pageref{onerealblock}. 

The converse statement  easily follows from direct verification (also one may notice that the above arguments are, in fact,  invertible and therefore $L$ and $g$ given by Theorem \ref{onerealblock} satisfy the compatibility equation \eqref{main} automatically). \qed



\section{Proof of Theorem \ref{onecomplexblock}: a pair of complex conjugate Jordan blocks} \label{proof3}


In this section we assume that $L$ has two complex conjugate eigenvalues $\lambda=a+ib$ and $\bar\lambda=a-ib$, $b\ne 0$ (each of geometric multiplicity one), so that  $L$ and $g$ can be simultaneously reduced to the following canonical forms
\begin{equation}
\label{canoncomplex}
L_{\mathrm{can}}=
\begin{pmatrix}
\begin{matrix} a & \!\!\! -b \\ b & a \end{matrix}   &  \begin{matrix} 1 & 0 \\ 0 & 1 \end{matrix} & & \\
  & \begin{matrix} a & \!\!\! -b \\ b & a \end{matrix} & \ddots  & \\
  &                  &    \ddots  &  \begin{matrix} 1 & 0 \\ 0 & 1 \end{matrix} \\
  &    &    & \begin{matrix} a & \!\!\! -b \\ b & a \end{matrix}
\end{pmatrix}
\quad
\mbox{and}
\quad
g_{\mathrm{can}} =
  \begin{pmatrix}
 &  & &   & 1 \\
 &  & & 1 &   \\
 &  & \iddots & & \\
 & 1 &      &  & \\
1 &  &      &  &
\end{pmatrix}
\end{equation}

By using the ``moving frame'' machinery as above, we can find the commutation relations between the elements of the canonical frame (associated with the canonical forms \eqref{canoncomplex} of $L$ and $g$) and describe the corresponding canonical coordinate system.  However, this approach leads to serious technical difficulties because the commutation relations turn out to be quite complicated.  To simplify them we will change the canonical forms of $L$ and $g$ in a certain way which is, in fact, motivated by  the splitting construction from \cite{splitting} which we recalled in \S\,\ref{sec:splitting}.
Namely we set, using Remark \ref{addremark}, page \pageref{addremark} with $P(t)=(t-\lambda)^n +(t-\bar\lambda)^n$:
\begin{equation}
\label{newcanon}
L_{\mathrm{can}}=L_{\mathrm{can}}^{\mathrm{old}}, \quad
g_{\mathrm{can}}=g_{\mathrm{can}}^{\mathrm{old}} \bigl((L_{\mathrm{can}}^{\mathrm{old}}-\lambda\cdot\Id)^n + (L_{\mathrm{can}}^{\mathrm{old}}-\bar \lambda\cdot\Id)^n\bigr),
\end{equation}
where $L_{\mathrm{can}}^{\mathrm{old}}$ and $g_{\mathrm{can}}^{\mathrm{old}}$ are as in \eqref{canoncomplex} and $n=\frac{1}{2}\dim M$.
Notice that the operator  $\bigl((L-\lambda\cdot\Id)^n + (L-\bar \lambda\cdot\Id)^n\bigr)$ is real, so $g_{\mathrm{can}}$ is a real symmetric matrix.

Let $e_1,f_1,e_2, f_2, \dots, e_n, f_n$ be the canonical frame associated with these (real) canonical forms \eqref{newcanon}. 
To simplify the commutation relations between them, we need one more modification. Namely, 
we pass from $e_i, f_i$ ($i=1,\dots, n$) to the natural complex frame
 $\xi_1, \dots, \xi_n, \eta_1, \dots, \eta_n$ by putting
\begin{equation}
\label{tocomp}
\xi_k=\frac{1}{2}(e_k-if_k) \quad \mbox{and}  \quad \eta_k = \frac{1}{2}(e_k + if_k)= \bar\xi_k.
\end{equation}

Thus, from now on we allow ourselves to use formal complex combinations of tangent vectors, i.e., we pass  
from the real tangent bundle $TM$ to its complexification $T^{\mathbb C} M$. In particular, we consider the complex vector fields $\xi = e + if$, $e,f\in \Gamma(TM)$,  and  treat them as differential operators on the space of complex-valued smooth functions $w(x)=u(x) + i v(x)$ on $M$:
$$
\xi (w) = (e(u)- f(v)) + i (e(v) + f(u)).
$$
The commutators of complex-valued vector fields and other objects of this kind are defined in the natural way.

According to Lemma \ref{uniqueness}, we now can uniquely reconstruct the commutation relations between the elements of the frame $\xi_1, \dots, \xi_n, \eta_1, \dots, \eta_n$ and information about derivatives of $\lambda$ and $\bar\lambda$ along these elements. Here is the result

\begin{proposition}
\label{aboutcomplex}Let $e_1,f_1,e_2, f_2, \dots, e_n, f_n$ be the canonical frame associated with canonical forms \eqref{newcanon}. Then the complex frame  
$$
\xi_1,\dots,\xi_n, \, \eta_1,\dots,\eta_n
$$
defined by \eqref{tocomp} satisfies the following properties:
\begin{enumerate}
\item $\xi_k$ and $\eta_m$ commute for all $k,m$;
\item $\xi_1, \dots, \xi_{n-1}$ commute and $\eta_1, \dots, \eta_{n-1}$ commute (in particular, all real vector fields $e_k$ and $f_m$ commute for all
$k,m\le n-1$);
\item the only nonzero derivatives are $\xi_n (\lambda)$ and $\eta_n(\bar \lambda)$;
\item nontrivial commutation relations are:
$$
\begin{array}{ll}
\, [\xi_1, \xi_n ] = -\ \xi_n (\lambda)\cdot  \xi_2, & [\eta_1, \eta_n]= -\  \eta_n(\bar \lambda)\cdot \eta_2, \\
\, [\xi_2, \xi_n ]= - 2\xi_n (\lambda) \cdot \xi_3, & [\eta_2, \eta_n]= -2 \eta_n(\bar \lambda)\cdot \eta_3, \\
\, [\xi_3, \xi_n]= - 3\xi_n (\lambda) \cdot \xi_4, & [\eta_3, \eta_n]= -3 \eta_n(\bar \lambda)\cdot \eta_4, \\
\dots & \dots \\
\, [\xi_{n-1}, \xi_n]= - (n-1)\xi_n (\lambda)\cdot \xi_n, & [\eta_{n-1}, \eta_n]= -(n-1) \eta_n(\bar \lambda)\cdot \eta_n.
\end{array}
$$
\end{enumerate}
\end{proposition}

Proof. One can find these relations
by straightforward (linear-algebraic) computation, but we shall give another proof based on the splitting construction (see \S \ref{sec:splitting})   and the  uniqueness lemma (Lemma \ref{uniqueness}, page \pageref{uniqueness}).

Before discussing the case of a complex Jordan block (more precisely, of two complex conjugate blocks), consider the case of two real Jordan blocks with distinct eigenvalues  as an illustrating example.
 Take two compatible pairs $(g_1, L_1)$ and $(g_2, L_2)$, each of which represents a single Jordan block with eigenvalue $\lambda_i \in \R$  (see the previous section for the complete description).  Let $\xi_1, \dots , \xi_n$ be the canonical frame for the first pair $(g_1, L_1)$ and $\eta_1,\dots, \eta_k$ for the second one $(g_2, L_2)$.  The gluing lemma (Theorem \ref{thm3}, page \pageref{thm3}) allows us to construct a new compatible pair $L$, $g$ by putting:
\begin{equation}
\label{gluing1}
L =
\begin{pmatrix}
L_1 & 0 \\ 0 & L_2
\end{pmatrix},
\quad
g=
\begin{pmatrix}
g_1 \chi_2 (L_1)& 0 \\ 0 & g_2 \chi_1(L_2)
\end{pmatrix}=
\begin{pmatrix}
g_1 & 0 \\ 0 & g_2
\end{pmatrix}
(\chi_1(L) + \chi_2(L))
\end{equation}
where $\chi_i(t)$ is the characteristic polynomial of $L_i$. The compatibility of $L_i$ and $g_i$, $i=1,2$, guarantees the compatibility of $L$ and $g$.

Now ask ourselves the converse question. Let $\xi_1, \dots , \xi_n, \eta_1,\dots, \eta_k$ be the canonical frame for a compatible pair $g$, $L$ having the (non-standard) canonical form \eqref{gluing1} with $L_i$, $g_i$ being the standard canonical forms as \eqref{canonicalJordan}.
What are the commutation relations between the elements of the frame and the conditions on the derivatives of the eigenvalues $\lambda_1$ and $\lambda_2$ along the frame?  We mean, of course, those relations which can be derived from the compatibility equation for $g$ and $L$.

Using the uniqueness result (Lemma \ref{uniqueness}), we immediately conclude that these relations will be exactly the same as for two separate Jordan blocks, namely,  $\xi_i$'s commute with $\eta_j$'s and the relations within each of these two groups will be those given in Lemmas \ref{firstcommutators} and \ref{secondcommutators} in \S  \ref{real}.

We now notice  that in this construction nothing changes, if we allow $\lambda_1$ and $\lambda_2$ to be complex conjugate, i.e., $\lambda_1=\lambda$ and $\lambda_2=\bar\lambda$ with  $\mathrm{Im}\lambda\ne 0$. The elements $\xi_1, \dots , \xi_n, \eta_1,\dots, \eta_k$ of the canonical frame will be, of course, vectors of the complexified tangent space $(T_PM)^{\mathbb C}$. The point is that Lemmas \ref{uniqueness}, \ref{en}, \ref{firstcommutators} and \ref{secondcommutators} are of purely algebraic nature and therefore can be applied for comlexified objects without any change.

If $e_1,f_1,e_2, f_2, \dots, e_n, f_n$ is the canonical frame associated with the (real) canonical forms \eqref{newcanon}, then in the complex frame
 $\xi_1, \dots, \xi_n, \eta_1, \dots, \eta_n$ defined by \eqref{tocomp}, 
 $L_{\mathrm{can}}$ and $g_{\mathrm{can}}$ take the form:
$$
L_{\mathrm{can}} \mapsto L'_{\mathrm{can}}  =
\begin{pmatrix}
\lambda  & 1         & & & & &  \\
          & \lambda  & \ddots & & &    \\
          &        &\ddots & 1 & & \\
  &        &         &\lambda &  & & \\
  &        &         &        &\bar\lambda &1 & \\
  & & &  &        & \bar\lambda & \ddots &  \\
& & & & & & \ddots & 1 \\
& & & & & & & \bar\lambda\\
\end{pmatrix} = \begin{pmatrix}  L_\lambda & 0 \\ 0 & L_{\bar\lambda} \end{pmatrix},
$$

$$
g_{\mathrm{can}} \mapsto g'_{\mathrm{can}} = \frac{1}{2}
\begin{pmatrix}
  & & \!\!\!\! -i & & & \\
  & \iddots & & & & \\
\!\!-i  & & & & & \\
  & & & & & i \\
  & & & &\iddots & \\
  & & & i & & \\
\end{pmatrix} \cdot \bigl(\chi_1 (L'_{\mathrm{can}})  + \chi_2 (L'_{\mathrm{can}})\bigr),
$$
where $\chi_1(t) = (t-\lambda)^n$ and $\chi_2(t) = (t-\bar\lambda)^n$ are the
characteristic polynomials of $L_\lambda$ and $L_{\bar\lambda}$ respectively.

According to Lemma \ref{uniqueness}, we can now uniquely reconstruct the commutation relations between the elements of the frame $\xi_1, \dots, \xi_n, \eta_1, \dots, \eta_n$ and information about the derivatives of $\lambda$ and $\bar\lambda$ along these elements. This reconstruction could be done by straightforward computation. Instead it suffices to notice that we are now essentially in the same situation as in the case of two real Jordan blocks,  but with $\lambda_1$ and $\lambda_2$ replaced by complex conjugate eigenvalues $\lambda$ and $\bar\lambda$, see formula \eqref{gluing1} and discussion around. So we can repeat the above argument to get the conclusion of Proposition \ref{aboutcomplex}. \qed

The canonical coordinate system  $x_1,y_1,x_2, y_2,\dots, x_n, y_n$ can now be reconstructed from these relations. To do it in the most natural way, we notice that $M$ carries a complex structure $J$ canonically associated with  $L$. Indeed, we can introduce $J$ in an invariant way as follows.  Let $L=L_{\mathrm s}+L_{\mathrm n}$ be the canonical decomposition of $L$ into semisimple and nilpotent parts. Then $J = \frac{1}{b} (L_{\mathrm s}-a)$ where $\lambda = a+ib$ is the eigenvalue of $L$. It is easy to see that $J^2 = -\Id$ and the integrability of $J$, i.e., vanishing of its Nijenhuis torsion $N_J$ follows from $N_L=0$.

\begin{Rem}
\label{aboutJ}
The fact that the Nijenhuis torsion $N_L$ of $L$ vanishes is well known (see e.g. \cite{benenti}), the implication
$N_L=0 \ \Rightarrow \ N_J=0$ can be verified directly.  Alternatively, one can use  \cite[Lemma 6]{splitting}. Indeed,  $J$ can be represented as the matrix function $J=f(L)$  that corresponds to the complex function $f: \mathbb C\setminus \R \to \mathbb C$ which is analytic, locally constant and defined  in the following way:  $f(a+ib)=i$ if  $b>0$ and $f(a+ib)=-i$ if  $b<0$. In a more general setting this construction of the canonical complex structure $J$ associated with $L$  is explained below in Lemma \ref{compstr}.
\end{Rem}

For the basis vectors $\xi_k$ and $\eta_k$ we have  $J\xi_k = i\xi_k$ and $J\eta_k = -i\eta_k$.  This means that for any holomorphic coordinate system
$z_1,\dots, z_n$,  the vectors $\xi_k$'s are linear combinations of $\partial_{z_k}$ and $\eta_k$'s are linear combinations of $\partial_{\bar z_k}$. Moreover, item 1 of Proposition \ref{aboutcomplex} says that the vector fields $\xi_1, \dots, \xi_n$ are holomorphic.

 From now on we can forget about $\eta_k$'s and work with $\xi_k$'s only. The following repeats the arguments for a real Jordan block but in the complex (holomorphic) setting. Since the holomorphic vector fields $\xi_1,\dots,\xi_{n-1}$ pairwise commute,
we can find a local complex coordinate system $z_1,\dots, z_n$  such that
$$
\xi_k = \partial_{z_k} \quad k=1,\dots, n-1
$$
Moreover, this coordinate system can be chosen in such a way that on the two-dimensional surface $z_1=z_2=\dots=z_{n-1}=0$, we have
\begin{equation}
\label{initialcomplex}
\xi_n = \partial_{z_n}.
\end{equation}

The eigenvalue $\lambda$ is a holomorphic function  (since $\xi_k$'s are holomorphic and $[\xi_{n-1}, \xi_n]= - (n-1)\xi_n (\lambda)\cdot \xi_n$). Moreover, $\lambda$ depends on $z_n$ only because $\xi_k(\lambda)=\partial_{z_k}(\lambda)=0$, $k=1,\dots,n-1$.

Now our goal is to determine the transition matrix between two bases $\xi_1, \dots, \xi_n$ and $\partial_{z_1}, \dots, \partial_{z_n},$.
As before, we consider the relation
$$
\partial_{z_n} = \sum_{k=1}^{n} a_{k-1} \xi_k =\sum_{k=1}^{n-1} a_{k-1} \partial_{z_k} + a_{n-1} \xi_n.
$$

Applying this differential operator to $\lambda$ we obtain:
$$
\lambda'_{z_n}=\frac{\partial \lambda}{\partial z_n} =  a_{n-1} \xi_n(\lambda)
$$

Next, computing $[\partial_{z_k}, \partial_{z_n}]$ by using the above formula, for $k=1,\dots, n-1$,  we see that $a_{k}$ depends on $z_{k}$ and $z_n$ only and satisfies:
$$
\frac{\partial a_{k}}{\partial z_{k}}= k a_{n-1}  \xi_n(\lambda) = k \, \lambda'_{z_n}.
$$

These equations can be easily solved:
$$
a_k=k z_{k} \, \lambda'_{z_n} + h_k (z_n, \bar z_n)
$$

And taking into account the initial conditions \eqref{initialcomplex}, we
conclude that
\begin{equation}
\label{complexa}
\aligned
& a_0=0,  \\
& a_1=  \lambda'_{z_n} z_1,  \\
& a_2=2  \lambda'_{z_n} z_2,  \\
& \dots                              \\
& a_{n-2}=(n-2) \lambda'_{z_n} z_{n-2}, \\
& a_{n-1}=1 + (n-1) \lambda'_{z_n} z_{n-1}.
\endaligned
\end{equation}

These formulas are, of course, identical to those for the real case. The only difference is that now we work with complex variables.

Now we have all the information to rewrite the formulas for $L$ and $g$ in the basis
$\partial_{z_i}$.  Notice that $L$ commutes with the complex structure $J$ and therefore $L$ can be treated as a complex operator.
The form $g$ is also compatible with $J$ in the sense that $g(Ju,v)=g(u,Jv)$ so that $g$ can be understood as the real part of the complex bilinear form on $T_p M$ treated as an $n$-dimensional complex space (w.r.t. $J$). This allows us to represent $L$ and $g$ by $n\times n$ complex matrices. In the canonical frame $\xi_1,\dots, \xi_n$ these matrices are:
$$
L^{\mathbb C}_{\mathrm{can}}  =
\begin{pmatrix}
\lambda  &  1       & &  \\
         & \lambda  & \ddots &     \\
         &          &\ddots  & 1  \\
  &        &         &\lambda  \\
\end{pmatrix}
$$

$$
 g^{\mathbb C}_{\mathrm{can}} =
\begin{pmatrix}
  & & \!\!\!\! -i  \\
  & \iddots &  \\
\!\!-i  & &  \\
\end{pmatrix} \cdot (L^{\mathbb C}_{\mathrm{can}}-\bar\lambda\cdot\Id)^n
$$

To determine $L$ and $g$  (more precisely their complex representations $L^{\mathbb C}$ and $g^{\mathbb C}$)  in the coordinates
$z_1,\dots, z_n$, we use the standard transformation:
$$
L^{\mathbb C}_{\mathrm{can}} \mapsto L^{\mathbb C}= C^{-1} L^{\mathbb C}_{\mathrm{can}} C, \quad
g^{\mathbb C}_{\mathrm{can}} \mapsto g^{\mathbb C}= C^{\top} g^{\mathbb C}_{\mathrm{can}} C
$$
with the transition matrix $C$
$$
(\partial_{z_1}, \dots, \partial_{z_n}) = (\xi_1,\dots , \xi_n) C,
\quad
C = \begin{pmatrix}
1 &  &   & &  a_0 \\
 & 1 &   &  & a_1 \\
& & \ddots & & \vdots \\
& &  & 1 & a_{n-2}\\
& & & & a_{n-1}
\end{pmatrix}
$$
where $a_k$ are defined by \eqref{complexa}.

Now a straightforward computation of $g^{\mathbb C}$ and $L^{\mathbb C}$ immeditely leads to the conclusion of Theorem \ref{onecomplexblock}, page \pageref{onecomplexblock}. \qed



\section{Applications: some global results} \label{proof4}


Here we give the proofs of Theorem  \ref{appl1}   and Corollary \ref{appl2}, page \pageref{appl1}.


\subsection{Proof of Theorem \ref{appl1}}


Consider two projectively equivalent pseudo-Riemannian metrics $g$ and $\bar g$ on $M$ and the $(1,1)$-tensor fields $L=L(g,\bar g)$ defined by \eqref{LG} which we repeat for convenience of the reader here:
$$
L^i_j= { \left|\frac{\det \bar g}{\det  g}\right|}^{\frac{1}{n+1}} \bar g^{ik} g_{kj}.
$$

Theorem \ref{appl1} can be reformulated as follows: 
\medskip

\centerline{{\it If $M$ is compact, then non real eigenvalues of $L$ are all constant.}}

The idea of the proof is very natural. As we know from \S\,\ref{proof3},  a complex non real  eigenvalue of $L$ is a holomorphic function in an appropriate coordinate system.  Roughly speaking, our proof is somehow equivalent to saying that ``a holomorphic function on a compact manifold has to be constant''.  However, to make sense out of this principle we have to deal with two issues:
\begin{itemize}
\item the complex structure $J$  (see Theorem \ref{onecomplexblock} and Section \ref{proof3}) is not globally defined;
\item the  eigenvalues of $L$ may collide and near the 
 collision points  they cannot be considered as well defined functions  (these points should be treated as branching points for eigenvalues).
\end{itemize}

To avoid these difficulties, we use the following two observations:
\begin{itemize}
\item a natural complex structure $J$ is well defined as soon as we have complex non real eigenvalues even at collision points,
\item the complex non real eigenvalues $\lambda_i$ of $L$  can be replaced by symmetric polynomials of them, like $\sum \lambda_i$, which are still holomorphic and well defined even at collision points.
\end{itemize}

These two ideas are formalized in the following lemma.
Let $\chi_L(t)$ be the characteristic polynomial of $L$. Clearly, the coefficients of $\chi_L(t)$ are smooth real functions on $M$.

\begin{Lemma}
\label{lemad1}
Let $\mu_0$ be a complex root of $\chi_L(t)$  of multiplicity $k$ at a point $p_0\in M$, $\mathrm{Im}\, \mu_0 > 0$.
Then in a neighborhood of $p_0$ there is a local coordinate system
$x_1,\dots, x_k, y_1,\dots, y_k, v_1, \dots, v_l$, $2k+l=\dim M$, such that the characteristic polynomial of $L$ admits the following factorization:
\begin{equation}
\label{factor}
\chi_L(t) = P_z(t)\cdot \bar P_z(t) \cdot Q_v(t)
\end{equation}
where
$$
P_z(t)=t^k + a_{k-1}(z)t^{k-1} + \dots + a_1(z) t + a_0(z)
$$
with coefficients $a_m(z)$ being holomorphic functions of the complex variables $z_j= x_j + i y_j$,
$$
\bar P_z(t)=t^k + \bar a_{k-1}(z)t^{k-1} + \dots + \bar a_1(z) t + \bar a_0(z),
$$
and
$$
Q_v(t)=t^l +  b_{l-1}(v)t^{l-1} + \dots +  b_1(v) t +  b_0(v)
$$
where $b_m(v)$ are smooth real valued functions of $v^1,\dots, v^l$,
and the polynomial $P_z(t)$ at the point $p_0$  takes the form $(t-\mu_0)^k$.
\end{Lemma}

{\it Proof.}
We first notice (as we did in our splitting construction \cite{splitting}) that in a neighborhood of $p_0$ the characteristic polynomial $\chi_L(t)$ can be uniquely factorized into two monic polynomials of degree $2k$ and $l$ respectively  with smooth real coefficients
$$
\chi_L(t)= \chi_1(t) \chi_2(t)
$$
in such a way that the roots of $\chi_1(t)$ at the point $p_0$ are $\mu_0$ and $\bar\mu_0$, both with multiplicity $k$. Locally, in a neighborhood of $p_0$ these polynomials do not have common roots. This factorization immediately lead (see  \cite[Theorem 2]{splitting}) to the existence of a coordinate system $u^1,\dots, u^{2k}$, $v^1, \dots, v^l$ in which $L$ splits into blocks each of which depends on its own group of coordinates:
$$
L= \begin{pmatrix}
L_1(u) & 0\\
0 & L_2(v)
\end{pmatrix},
$$
so that $\chi_1(t)$ and $\chi_2(t)$ are the characteristic polynomials of $L_1$ and $L_2$. In other words, locally we may think of $M$ with $L$ as a direct product  $(M_1,L_1)\times (M_2,L_2)$ of two ``independent'' manifolds with $(1,1)$-tensor fields on them.  We put $Q_v(t)=\chi_2(t)$ and continue working,  from now on, with the first factor $(M_1, L_1)$ only.

In a neighborhood of $p_0$, the characteristic polynomial $\chi_1(t)$ of $L_1$ admits a further factorization:
$$
\chi_1(t) = P_u(t)\bar P_u(t)
$$
into two complex conjugate polynomials with smooth complex valued coefficients satisfying the required property: at the point $p_0$ we have
$P_u(t)=(t-\mu_0)^k$.  So far this construction is purely algebraic. But now we need to pass from real coordinates $u$ to complex coordinates
$z = x + i y$ in such a way that the coefficients of $P_u(t)$ become holomorphic functions of $z$.

First we construct a complex structure $J$ on $M_1$ canonically associated with $L_1$.   Since this construction seems to be quite natural and perhaps useful in other geometric problems, we explain it in a more general context,  referring to \cite{splitting}, Section 1.3 for technical details.

Consider a locally constant complex function $f: \mathbb C\setminus \R \to \mathbb C$  defined by:
$$
f(z) = \begin{cases}
\ \ i, \quad \mbox{if } \ \mathrm{Im}\,  z >0,\\
-i, \quad \mbox{if }  \ \mathrm{Im}\,z <0.
\end{cases}
$$
Since $f$ can be uniformly approximated (with all derivatives up to any fixed order) by polynomials with real coefficients on every compact subset $K\subset \mathbb C\setminus\R$,  we are also allowed to consider $f$ as a real analytic {\it matrix} function $f: U \to \mathrm{End}(\R^{n})$, where $U\subset  \mathrm{End}(\R^n)$ is a subset of all operators with no real eigenvalues (in particular, here $n=2k$).  

Now, let $L$ be a smooth $(1,1)$-tensor field on a smooth manifold $M$    with no real eigenvalues (like $L_1$ on $M_1$ in our case),  i.e., the spectrum of $L$ at every point $x\in M$ belongs to $\mathbb C \setminus\R$. Then we can build a new $(1,1)$-tensor field $J$ on $M$ by setting $J= f(L)$ pointwise. 

  \begin{Lemma} 
  \label{compstr}
  The $(1,1)$-tensor field $J$ satisfies the following properties:
  \begin{enumerate}
  \item $J$ is smooth;
  \item $J^2 = -\mathbf{1}$, i.e. $J$ is an almost complex structure on $M$;
  \item $JL=LJ$, i.e. $L$ is a complex linear operator w.r.t. $J$;
  \item if the Nijenhuis torsion $N_L$ of $L$ vanishes, then  $J$  is integrable and hence is a complex structure on $M$. 
 \end{enumerate}
\end{Lemma}

 {\it Proof.} The smoothness of $J$ follows from the fact that  $f: U \to \mathrm{End}(\R^{n})$ is a real analytic matrix function and $L$ depends smoothly of $x\in M$. Next, items 2 and 3 are purely algebraic.  Indeed,  the scalar identity $f^2(z)\equiv -1$ implies the matrix identity $J^2 = f(L) f(L) = - \mathbf{1}$, and $L$ commutes with $f(L)$  for any matrix function $f$. Finally, the integrability of $J$, i.e.,  the fact that $N_J\equiv 0$, is a particular case of Lemma~6  from \cite{splitting}. \qed

\begin{Rem}\label{newrem}
An equivalent definition of $J$, as a function of $L$, is as follows. Let $L: V^{2k} \to V^{2k}$ be a real linear operator with no real eigenvalues.  We consider the decomposition of $V^{\mathbb C}$ into two $L$-invariant subspaces $V^+ \oplus V^-$ corresponding to the eigenvalues of $L$ with positive and negative imaginary parts respectively. Such a decomposition is obviously unique.  Now we define $J$ to be the multiplication by $i$ on $V^{+}$ and multiplication by $-i$ on $V^-$. It is easy to see that $V$ as a subspace of $V^{\mathbb C}$ is $J$-invariant, i.e.,  $J$ gives a well-defined operator on $V$, satisfying $J^2=-\Id$ and commuting with $L$.
In particular, in the case of a single complex Jordan block, the complex structure $J$ canonically associated with $L$ coincides with the one we used in \S\,\ref{proof3} (cf. Remark \ref{aboutJ}, page \pageref{aboutJ}). Moreover,  in a neighborhood of a regular point, $J$  is the direct sum of the complex structures constructed for each individual $(\lambda_i,\bar\lambda_i)$-block by the method explained in \S\,\ref{proof3}.  
 \end{Rem}

Let us come back to the proof of Lemma \ref{lemad1}.  By applying Lemma \ref{compstr} to the $(1,1)$-tensor field $L_1$ on $M_1$, we construct the complex structure $J$ on $M_1$ canonically associated with $L_1$. 
Next, we need to show that the coefficients of $P_z(t)$  (obtained from $P_u(t)$ by replacing the coordinates $u$ with $z$) are holomorphic with respect to $J$. Though it can be done independently, we shall easily derive this property from \S\,\ref{proof3}.

Indeed, in  \S\,\ref{proof3} we have shown that in a neighborhood of every regular  point each complex non real eigenvalue $\lambda_i$ of $L$ is a holomorphic function (in the ``singular'' case studied in \S\,\ref{affine}, $\lambda_i$ is constant, so this property holds automatically)  w.r.t. the complex structure associated with the $(\lambda_i,\bar\lambda_i)$-block. Taking into account Remark \ref{newrem}, we see that each $\lambda_i$ is holomorphic w.r.t. $J$.
On the other hand the coefficients $a_j(z)$ of $P_z(t)$ are symmetric polynomials in $\lambda_i$,  so they are holomorphic at each regular  point too. Now it remains to notice that regular  points form an open dense subset and $a_j(z)$ are smooth everywhere. This obviously implies that
$a_j(z)$ are holomorphic on the whole neighborhood of $p_0$.  This completes the proof. \qed

We also shall use the following almost obvious statement.

\begin{Lemma}
\label{lemad2}
Let $P_{z}(t) = t^k + a_{k-1}(z)t^{k-1} + \dots + a_1(z) t + a_0(z)$ be a polynomial in $t$ whose coefficients are holomorphic functions of $z=(z_1,\dots, z_k)\in U$, where $U\subset \mathbb C^k$ is an open connected domain.

Assume that at some point $z_0\in U$ the polynomial takes the form $P_{z_0}(t) = (t-\mu_0)^k$
and at any other point $z\in U$ all the roots $\lambda_i(z)$ of $P_z(t)$ satisfy the condition
$$
\mathrm{Im}\, \lambda_i (z) \le  c=\mathrm{Im}\, \mu_0, \quad  i=1,\dots, k.
$$
Then $\lambda_i(z)\equiv \mu_0$ for all $z\in U$, i.e., the roots of $P_z(t)$ are all constant and equal to $\mu_0$. In particular,
$P_z(t)\equiv (t-\mu_0)^k$ on $U$.
\end{Lemma}

{\it Proof.}
Consider the sum $\sum_{i=1}^k \lambda_i(z)$ of the roots of $P_z(t)$.
Since $\sum_{i=1}^k \lambda_i(z)=-a_{k-1}(z)$, this sum is a holomorphic function on $U$.
On the other hand, we see that $\mathrm{Im} (- a_{k-1}(z)) = \mathrm{Im}\,\left( \sum_{i=1}^k \lambda_i(z) \right) \le k\cdot c$
and $\mathrm{Im} (- a_{k-1}(z_0))=k\cdot c$, i.e., the imaginary part of the holomorphic function $- a_{k-1}(z)$ attains a maximum at a certain point $z_0\in U$.  This implies (by the maximum principle) that $a_{k-1}(z)$ is constant on $U$. From this, in turn, it is easy to derive that the imaginary part of each $\lambda_i(z)$ and, therefore, $\lambda_i(z)$ itself is constant.

We are now ready to complete the proof of Theorem \ref{appl1}.

Consider the roots $\lambda_1(p), \dots, \lambda_n(p)$, $n=\dim M$ of the
characteristic polynomial $\chi_L(t)$ at $p\in M$ and let
$$
c=\max_{p\in M, \ i=1,\dots, n} \mathrm{Im}\, \lambda_i(p)
$$

We assume that some complex eigenvalues exist, so $c>0$.

Since the roots $\lambda_i(p)$ depend on $p$ continuously (in a natural sense) and $M$ is compact, then $c$ is attained, i.e., there is a point  at which  $\chi_{L}(t)$ has a complex root $\mu_0$ such that $\mathrm{Im}\,\mu_0=c$. In general, $\mu_0$ may have different multiplicities at different points. Let $k$ be maximal multiplicity of $\mu_0$ on $M$.

Consider the following subset $A \subset M$:
$$
A=\{ q\in M~|~ \mbox {$\mu_0$ is a root of  $\chi_{L}(t)$ of multiplicity $k$ at the point $q$}\}.
$$

By our assumption, $A$ is nonempty and 
as the multiplicity of $\mu_0$ is upper semi-continuous and  $k$ is its maximum, $A$ is closed.
Let us show that $A$ is open. Indeed, let $p_0\in A$. We first apply Lemma \ref{lemad1} at this point to get the factorization \eqref{factor} in some neighborhood $U(p_0)$ and then apply
Lemma \ref{lemad2} to see that $\chi_L(t) = (t-\mu_0)^k (t-\bar\mu_0)^k Q_v(t)$ on  $U(p_0)$. In other words, $\mu_0$ is a root of $\chi_L(t)$ of multiplicity $k$ for all points $p\in U(p_0)$, i.e. $U(p_0)\subset A$ and therefore $A$ is open.

Thus, $A$ is open, closed and non-empty. Hence, $A=M$ and we see that $\chi_L(t)=(t-\mu_0)^k(t-\bar\mu_0)^k Q_v(t)$ everywhere on $M$. In other words, $\mu_0$ is a constant complex eigenvalue of $L$ of multiplicity $k$ on the whole manifold $M$.

If  $Q_v(t)$  has some other non real roots at some points of $M$, we simply repeat the same argument to show that these roots have to be constant. This completes the proof of  the first statement 
Theorem \ref{appl1}. 

In order to proof that the multiplicity  of a  non real eigenvalue is the same at every point, it is sufficiently to observe that otherwise there must  be a non real nonconstant eigenvalue of $L$, which contradicts the proven part of Theorem \ref{appl1}. \qed



\subsection{Proof of Corollary \ref{appl2}}


We need to prove the following result:

{\it 
Let $M^3$ be a closed connected 3-dimensional manifold. Suppose  $g$ and $\bar g$ are geodesically equivalent metrics on it and $L$, given by \eqref{L}, has a non real eigenvalue  at least  at
one point.
Then, $M^3$ can be finitely covered by the 3-torus.}

 Without loss of generality we assume that the metric $g$ has signature $(-,+,+)$.
 Then, by Theorem \ref{appl1},  page \pageref{appl1}, a non real eigenvalue of $L$  is a constant; we denote it  by $\alpha+ i\beta$.  The complex conjugate number $\alpha- i\beta$ is also an  eigenvalue of $L$;
 the remaining third eigenvalue will be denoted by $\lambda$; it is a (smooth) real valued function on the manifold.

At every point $p\in M$, let us consider a basis $\{v_1,v_2, v_3\}$
in  $T_pM$  such that  in this basis the matrices of $g$ and $L$ are given by
\begin{equation} \label{g3}
g = \begin{pmatrix}(\lambda - \alpha)^2 + \beta^2  & & \\ & -\beta & \alpha- \lambda \\ &\alpha-\lambda & \beta   \end{pmatrix}, \ \ L   = \begin{pmatrix}\lambda & & \\ & \alpha & \beta  \\ &-\beta & \alpha  \end{pmatrix}.
\end{equation}

The existence of such a basis follows from \cite[Theorem 12.2]{lancaster}; it is an easy exercise to show that  the basis  is unique up to the transformations $v_1\mapsto - v_1$; $(v_2, v_3)\mapsto (-v_2, -v_3)$.

Now, consider the positive definite Euclidean structure at $T_pM$ such that this basis is orthonormal. This Euclidean structure does not depend on the freedom in the choice of the basis and is therefore well defined.   It smoothly depends on the point $p$  and therefore generates a Riemannian metric on $M$, which we denote by $g_0$. Let us show that the metric $g_0$ is flat.

In order to do it, we will use our description of compatible  pairs $(g,L)$. As we explained  in \S \ref{sec:splitting},  in a neighborhood of  every point the metric $g$ could be obtained by gluing $(I, h_1, L_1)$  and $(U^2, h_2, L_2)$
where \begin{itemize}\item
$I$ is one-dimensional,  the metric $h_1$ is positive definite and the eigenvalue of $L_1$ is $\lambda$.

\item  $U^2$ is two-dimensional, $h_2$ has signature $(-,+)$,  the eigenvalues of  $L_2$ are $\alpha+ i\beta, \alpha-i\beta$ and $h_2$ and $L_2$ are compatible.
\end{itemize}

Then,  for a  certain choice of the coordinate $x_1$  on $I$
 the metric $h_1$  is $(dx_1)^2$ and  the only component of the ($1\times 1$)-matrix of  $L$ is $\lambda(x_1)$.
 Now, since $L_2$ is compatible with $h_2$ and since the trace of $L_2$ is  constant, we conclude from \eqref{main} that $L_2$ is covariantly constant with respect to $h_2$. Then, $h_2$ is obviously flat and in a certain (local)
 coordinate system $(x_2,x_3)$ on $U^2$ the metric $h_2$ and the $(1,1)$-tensor $L$ are given by the matrices
$$
h_2 = \begin{pmatrix}0 & 1 \\ 1& 0  \end{pmatrix}, \ \  L_2= \begin{pmatrix}  \alpha & \beta  \\ -\beta & \alpha  \end{pmatrix}.
$$
  Applying the gluing construction to $(I, h_1, L_1)$  and $(U^2, h_2, L_2)$, we obtain that
the metric $g$ and the tensor $L$ are given by
$$
g= \begin{pmatrix}(\lambda(x_1) - \alpha)^2 + \beta^2& & \\ & -\beta & \alpha - \lambda(x_1)\\ & \alpha - \lambda(x_1)& \beta \end{pmatrix}, \ \ L= \begin{pmatrix}  \lambda(x_1) & & \\ &\alpha & \beta  \\ &-\beta & \alpha  \end{pmatrix}.$$
We see that the vector fields
$v_1= \tfrac{\partial }{\partial x_1}$, $v_2= \tfrac{\partial }{\partial x_2}$, $v_3= \tfrac{\partial }{\partial x_3}$ form a basic such that $g$ and $L$ are as in \eqref{g3} implying that the metric $g_0$ is given by $g_0= (dx_1)^2+ (dx_2)^2+(dx_3)^2$ and therefore is flat.

Thus, there exists a flat Riemannian  metric on $M^3$. Then, the manifold is a  3-dimensional Bieberbach manifold (i.e., is a quotient of $\mathbb{R}^3$ modulo a freely acting crystallographic group) and  can be finitely covered by the torus $T^3$ as we claimed. \qed

\weg{
\subsubsection*{Acknowledgements}
We are very grateful to the anonymous referee for many helpful comments and suggestions.
We thank DAAD (Programm Ostpartnerschaft) and DFG (GK 1523) for partial financial support. }

 \section*{  Appendix. What is a  function of a matrix or of a (1,1)-tensor (following \cite{splitting,higham})?}

In this appendix we give more details on the construction of the complex structure $J$ (canonically associated with $L$) that was used in the proofs of Theorems   \ref{onecomplexblock}, \ref{appl1}. The arguments below partially repeat those in the main part  of this paper and the paper is, in fact,  self-contained without this appendix. However from  discussions with colleagues and from remarks of the anonymous referee,  we understood that this construction looks slightly strange at first glance, so we 
repeat it separately here with detailed explanation.

Let $L$ be a $(1,1)$-tensor field with vanishing Nijenhuis torsion,  i.e. $N_L \equiv 0$. 

Step 1.  Let $P(t)$ be an arbitrary polynomial (with constant real coefficients).  Then $N_L \equiv 0$ implies $N_{P(L)} \equiv 0$, see \cite[first part of Lemma 6]{splitting}.   

Step 2. Instead of a polynomial $P(t)$, one can take any real  analytic entire  function $f$ like $\exp$, $\sin$, $\cos$, etc.  In this case, $f(L)$ is well-defined as the power series,    and  
$N_L \equiv 0$ implies $N_{f(L)} \equiv 0$. The proof is evident. 

Step 3. Thus,  we come to the natural question about the class of analytic functions $f$  which can be used in this context, i.e. for which, in particular,   the expression $f(L)$ makes sense.  Of course, every entire analytic  function $f: \mathbb C \to \mathbb C$  with the property $f(\bar z)= \overline{f(z)}$ (which is, for entire functions, equivalent to real-analyticity) 
is suitable for this purpose, but unfortunately this class of functions is too small for our goals.  Notice that although the operators $L$ are real, the functions must be defined on the complex plane, as $L$ might have complex eigenvalues.

Step 4.  The most reasonable class of functions, satisfying all the properties that we need,  is described by Mergelyan's theorem, see e.g. \cite{mergelyan}.  The below definition of such functions is copied from \cite[\S 1.3]{splitting}.

Suppose a compact set   $K \subseteq \mathbb{C} $,  a function $f: {K} \to \mathbb{C} $ and  an open subset $U\subset \mathrm{End}(\R^n) \simeq \mathrm{gl}(n,\R)$ satisfy the following conditions:

\begin{enumerate} 

\item $\mathbb{C}\setminus K$ is connected (we do not require that $K$  is connected); \label{i}

 \item  $f:K\to \mathbb{C}$ is  a continuous function and  the restriction $f |_{\mathrm{Int}(K)}$ is holomorphic; \label{ii}

\item $K$ is symmetric with respect to the $x-$axes: for every  $z\in K$ its conjugate $\bar z $ also lies in $K$; \label{iii}

\item for every $z\in K$,  $f(z)= \bar f(\bar z)$, where the bar ``$\bar{ \hspace{1ex} }$''    denotes the complex conjugation;  \label{iiii}

\item  for  every $L\in U$  we have $\mathrm{Spectrum}(L)\subset \mathrm{Int}(K)$. \label{iiiii} 
  
  \end{enumerate}

Step 5.    Under the above assumptions (\ref{i}--\ref{iiiii}), one can naturally define  a matrix function $$f: U \to \mathrm{End}(\R^n)$$ by using Mergelyan's theorem which states that  the function $f: K\to \mathbb C$  can be uniformly approximated by {\it real} polynomials $P_i$  (i.e., polynomials of complex variable $z$, but with real coefficients).  The integral Cauchy formula immediately implies that the sequence $P_i$ converges to $f$ together with all derivatives (and  moreover  uniformly on every  disc $B\subset \mathrm{Int} K$).

 We define $f(L)= \lim_{i\to \infty }P_i(L)$. It is not hard, see  
  for example   \cite[\S\S 1.2.2 -- 1.2.4]{higham},  to prove the following
  
  \begin{proposition}
  1) The limit  $f(L)= \lim_{i\to \infty }P_i(L)$ exists (and does not depend on the choice of the sequence $P_i$) for every $L\subset U$, where  $U\subset \mathrm{End}(\R^n)$ satisfies (\ref{iiiii}).

  2) The function $f: U \to \mathrm{End}(\R^n)$ so obtained is real analytic   (in the usual sense, i.e., as a map between  vector spaces).
  
  3) The partial derivatives  also converge, i.e.  $\dfrac{\partial f}{ \partial l^\alpha_\beta} (L) =  \lim_{i\to \infty }\dfrac{\partial p_i}{\partial l^\alpha_\beta}(L)$,   where $l^\alpha_\beta$ denote entries of $L$.  The same holds for all derivatives of higher order and this convergence in uniform on every compact subset $V\subset U$.
  
  4) $\mathrm{Spectrum}(f(L))= f( \mathrm{Spectrum}(L))$.
  \end{proposition}

Step 6.   This construction can be naturally extended to $(1,1)$-tensor fields. Namely, if $L=L(x)$ is a smooth $(1,1)$-tensor field on a smooth manifold $M$, then we can construct a new $(1,1)$-tensor field $f(L)$ on $M$ 
    (which will smoothly depend on $x$  being a composition of a smooth ``map'' $L: M\to \mathrm{End}(TM)$ and an analytic map $f$).

  The statement from \cite{splitting}  which we use in the present paper is
  \begin{proposition}[follows from Lemma 6 of \cite{splitting}] \label{fL} 
  Let $L$ be a smooth $(1,1)$-tensor field on $M$  and $f$ be a matrix function (of the above type) such that the $(1,1)$-tensor field $f(L)$ is well defined on $M$.  Then  $N_L=0$ implies $N_{f(L)}=0$.
  \end{proposition}  
  
  The proof is obvious. We can approximate $f(L)$  by polynomials $P_i(L)$  in such a way that all partial derivatives  $\partial_{x_k} P_i(L)$  also converge to $\partial_{x_k} f(L)$. Since the Nijenhuis torsion is an algebraic expression involving the components of a $(1,1)$-tensor field and the partial derivatives of them,  the statement follows.

  In \cite{splitting}  we used this statement to prove that projectors on some invariant subspaces of $L$ have zero Nijenhuis torsion.  In the present paper,  we use this fact  to construct a canonical complex structure on $M$ in the following situation.
  
  Let $L$ be a $(1,1)$-tensor field on $M$ with no real eigenvalues.  Then we consider  the following function on $\mathbb C \setminus \R$:
  $$
  f(z) = \left\{  \begin{array}{rl} 
  i \ , & \mbox{if } \mathrm{Im}\, z >0, \\
  - i \ , & \mbox{if } \mathrm{Im}\, z <0. 
  \end{array}
  \right.
  $$
 Formally,  we need to define this function on a compact subset $K \subset  \mathbb C \setminus \R$ satisfying the assumptions (1), (3) and (5) above, but this $K$ can be arbitrary large and the result does not depend on the choice of $K$ so the construction works on the whole set $\mathbb C \setminus \R$. 
  
  This function is locally constant (but not constant as the set $\mathbb C \setminus \R$ consists of two connected components and the values of the function on different connected components are different)
   and analytic on $\mathbb C \setminus \R$. It is easy to see that all properties \ref{i}--\ref{iiiii} are fulfilled  (formally speaking for any symmetric compact subset $K \subset \mathbb C\setminus \R)$) and therefore  the function $f(L)$ is defined for every $L$ with no real eigenvalues.
  
  It is easy to see that the operator $f(L)$ is a complex structure  (indeed  $f^2 = -1$ which implies that $ f(L) f(L) = f^2(L) =-\mathbf{1}$).
  
In our paper, the  function $f$  comes into this construction in three different ways,  namely
  \begin{enumerate}
  \item $f$ is a complex function defined on a subset $K\subset \mathbb C$ of a complex plane, i.e., $f: K\to \mathbb C$.  The function $f: \mathbb C\setminus \R\to \mathbb C$ that we use in our paper is complex analytic (holomorphic on $\mathbb C\setminus \R$) and locally constant.
   
  \item  We can also consider $f$ as a matrix function $f: U \to \mathrm{End} (\R^n)$, where $U\subset \mathrm{End} (\R^n)$ is a subset  satisfying (\ref{iiiii}). This function is real,  analytic and non constant ($J=f(L)$ depends on $L$ in a very nontrivial way, it is not even a polynomial!).
  
  \item Finally, we can consider the expression $f(L)$ for a smooth $(1,1)$-tensor field $L$ on $M$ and then $f(L)$  as a ``function'' on a manifolds $M$  (section of  $\mathrm{End}(TM)$). From this point of wiev, $f(L)$ is a smooth function/section/tensor field in the sense of smooth dependence on local coordinates $x_1,\dots, x_n$ on $M$.

  \end{enumerate}

  Now if $L$ is a $(1,1)$-tensor field on $M$ with no real eigenvalues at any point $x\in M$, then $J=f(L)$ is an almost complex structure on $M$ (canonically associated with $L$) and  Proposition \ref{fL} immediately implies
  
\begin{Cor}
If $N_L =0$, then $N_J=N_{f(L)}=0$ and therefore $J$ is a complex structure on $M$.
\end{Cor}


\end{document}